\documentclass[11pt]{article}
\usepackage{}
\usepackage{mathrsfs}
\usepackage{amsfonts}
\usepackage{amsfonts,amssymb,mathrsfs}
\usepackage{amsmath,amscd}
\usepackage[dvips]{graphicx}
\usepackage[all]{xy}
\usepackage{graphicx}
\usepackage{fancyhdr}
\usepackage{mathptm,pslatex}
\usepackage{amsthm}

\begin{document}

\sloppy
\newtheorem{Def}{Definition}[section]
\newtheorem{Bsp}{Example}[section]
\newtheorem{Prop}[Def]{Proposition}
\newtheorem{Theo}[Def]{Theorem}
\newtheorem{Lem}[Def]{Lemma}
\newtheorem{Koro}[Def]{Corollary}
\theoremstyle{definition}
\newtheorem{Rem}[Def]{Remark}

\newcommand{\add}{{\rm add}}
\newcommand{\gd}{{\rm gl.dim }}
\newcommand{\dm}{{\rm dom.dim }}
\newcommand{\E}{{\rm E}}
\newcommand{\Mor}{{\rm Morph}}
\newcommand{\End}{{\rm End}}
\newcommand{\ind}{{\rm ind}}
\newcommand{\rsd}{{\rm res.dim}}
\newcommand{\rd} {{\rm rep.dim}}
\newcommand{\ol}{\overline}
\newcommand{\rad}{{\rm rad}}
\newcommand{\soc}{{\rm soc}}
\renewcommand{\top}{{\rm top}}
\newcommand{\pd}{{\rm proj.dim}}
\newcommand{\id}{{\rm inj.dim}}
\newcommand{\Fac}{{\rm Fac}}
\newcommand{\fd} {{\rm fin.dim }}
\newcommand{\DTr}{{\rm DTr}}
\newcommand{\cpx}[1]{#1^{\bullet}}
\newcommand{\D}[1]{{\mathscr D}(#1)}
\newcommand{\Dz}[1]{{\mathscr D}^+(#1)}
\newcommand{\Df}[1]{{\mathscr D}^-(#1)}
\newcommand{\Db}[1]{{\mathscr D}^b(#1)}
\newcommand{\C}[1]{{\mathscr C}(#1)}
\newcommand{\Cz}[1]{{\mathscr C}^+(#1)}
\newcommand{\Cf}[1]{{\mathscr C}^-(#1)}
\newcommand{\Cb}[1]{{\mathscr C}^b(#1)}
\newcommand{\K}[1]{{\mathscr K}(#1)}
\newcommand{\Kz}[1]{{\mathscr K}^+(#1)}
\newcommand{\Kf}[1]{{\mathscr  K}^-(#1)}
\newcommand{\Kb}[1]{{\mathscr K}^b(#1)}
%\stackrel{\sim}
\newcommand{\modcat}{\ensuremath{\mbox{{\rm -mod}}}}
\newcommand{\Modcat}{\ensuremath{\mbox{{\rm -Mod}}}}
\newcommand{\stmodcat}[1]{#1\mbox{{\rm -{\underline{mod}}}}}
\newcommand{\pmodcat}[1]{#1\mbox{{\rm -proj}}}
\newcommand{\imodcat}[1]{#1\mbox{{\rm -inj}}}
\newcommand{\opp}{^{\rm op}}
\newcommand{\otimesL}{\otimes^{\rm\bf L}}
\newcommand{\rHom}{{\rm\bf R}{\rm Hom}}
\newcommand{\projdim}{\pd}
\newcommand{\Hom}{{\rm Hom}}
\newcommand{\Coker}{{\rm coker}\,\,}
\newcommand{ \Ker  }{{\rm Ker}\,\,}
\newcommand{ \Img  }{{\rm Im}\,\,}
\newcommand{\Ext}{{\rm Ext}}
\newcommand{\StHom}{{\rm \underline{Hom} \, }}

\newcommand{\gm}{{\rm _{\Gamma_M}}}
\newcommand{\gmr}{{\rm _{\Gamma_M^R}}}

\def\vez{\varepsilon}\def\bz{\bigoplus}  \def\sz {\oplus}
\def\epa{\xrightarrow} \def\inja{\hookrightarrow}

\newcommand{\lra}{\longrightarrow}
\newcommand{\lraf}[1]{\stackrel{#1}{\lra}}
\newcommand{\ra}{\rightarrow}
\newcommand{\dk}{{\rm dim_{_{k}}}}

{\Large \bf
\begin{center}
Inductions and restrictions for stable equivalences of Morita type
\end{center}}
\medskip

\centerline{{\bf Hongxing Chen, Shengyong Pan} and {\bf Changchang
Xi$^*$}}
\begin{center} School of Mathematical Sciences, Beijing Normal University, \\
Laboratory of Mathematics and Complex Systems, \\
100875 Beijing, People's Republic of  China \\ E-mail:  chx19830818@163.com \quad panshy1979@mail.bnu.edu.cn \quad xicc@bnu.edu.cn\\
\end{center}
\bigskip

\renewcommand{\thefootnote}{\alph{footnote}}
\setcounter{footnote}{-1} \footnote{ $^*$ Corresponding author.
Email: xicc@bnu.edu.cn; Fax: 0086 10 58802136; Tel.: 0086 10
58808877.}
\renewcommand{\thefootnote}{\alph{footnote}}
\setcounter{footnote}{-1} 
\renewcommand{\thefootnote}{\alph{footnote}}
\setcounter{footnote}{-1} 

%\date{}

\begin{abstract}
In this paper, we present two methods, induction and restriction
procedures, to construct new stable equivalences of Morita type.
Suppose that a stable equivalence of Morita type between two
algebras $A$ and $B$ is defined by a $B$-$A$-bimodule $N$. Then, for
any finite admissible set $\Phi$ and any generator $X$ of the
$A$-module category, the $\Phi$-Auslander-Yoneda algebras of $X$ and
$N\otimes_AX$ are stably equivalent of Morita type. Moreover, under
certain conditions, we transfer stable equivalences of Morita type
between $A$ and $B$ to ones between $eAe$ and $fBf$, where $e$ and
$f$ are idempotent elements in $A$ and $B$, respectively.
Consequently, for self-injective algebras $A$ and $B$ over a field
without semisimple direct summands, and for any $A$-module $X$ and
$B$-module $Y$, if the $\Phi$-Auslander-Yoneda algebras of $A\oplus
X$ and $B\oplus Y$ are stably equivalent of Morita type for one
finite admissible set $\Phi$, then so are the
$\Psi$-Auslander-Yoneda algebras of $A\oplus X$ and $B\oplus Y$ for
{\it every} finite admissible set $\Psi$. Moreover, two
representation-finite algebras over a field without semisimple
direct summands are stably equivalent of Morita type if and only if
so are their Auslander algebras. As another consequence, we
construct an infinite family of algebras of the same dimension and
the same dominant dimension such that they are pairwise derived
equivalent, but not stably equivalent of Morita type. This answers a
question by Thorsten Holm.
\end{abstract}

\tableofcontents

\section{Introduction}

In the representation theory of algebras and groups, there are three
fundamental equivalences: Morita, derived and stable equivalences.
Roughly speaking, the first two are induced from tensor products of
bimodules or two-sided complexes, thus there is a corresponding
Morita theory for each (see \cite{morita, Rickard, keller}), while
the last one seems not yet to be well understood in this way, and
therefore a Morita theory for stable equivalences is missing.
Recently, a special class of stable equivalences, called stable
equivalences of Morita type, are introduced by Brou\'e in modular
representations of finite groups. They are induced by bimodules,
have features of a Morita theory, and are shown to be of great
interest in modern representation theory since they preserve many
homological and structural invariants of algebras and modules (see,
for example, \cite{B, DM, krause, L, P, C3, C2}). In order to
understand this kind of equivalences, one has to know, first of all,
examples and basic properties of stable equivalences of Morita type
as many as possible. So, one of crucial questions in the course of
studying these equivalences is:

{\bf Question:} How to construct stable equivalences of Morita type
for finite-dimensional algebras ?

\smallskip
\noindent Up to date, only a few  methods using trivial extensions,
one-point extensions and endomorphism algebras have been known in
\cite{R3, LX1, LX2, LX3}. Of course, Rickard's result that the
existence of derived equivalences for self-injective algebras
implies the one of stable equivalences of Morita type provides
another way to construct stable equivalences of Morita type. This
method, however, is no longer true for general finite-dimensional
algebras (see \cite{hx3} for some new advances in this direction).
So, a systematical method for constructing stable equivalences of
Morita type seems not yet to be available.

In this paper, we shall look for a more general and systematical
answer to this question, and present two methods, called induction
and restriction procedures, to construct new stable equivalences of
Morita type for general finite-dimensional algebras. Here our
induction procedure has two flexibilities, one is the choice of
generators, and the other is the one of finite admissible sets. Thus
this construction provides a large variety of stable equivalences of
Morita type.

To state our first main result, let us recall the definition of
$\Phi$-Auslander-Yoneda algebras in \cite{HAU}. Let $A$ be a
finite-dimensional algebra and $X$ an $A$-module. Then, for an
admissible set $\Phi$ of natural numbers, there is defined an
algebra $\E^{\Phi}_A(X)$, called the $\Phi$-Auslander-Yoneda algebra
of $X$ in \cite{HAU}, which is equal to $\bigoplus_{i\in
\Phi}\mbox{Ext}^j_A(X,X)$ as a vector space, and its multiplication
is defined in a natural way (see Subsection \ref{adm} below for
details). Our main result for inductions reads as follows:

\begin{Theo}{\rm (The Induction Procedure)}

Suppose that $A$ and $B$ are finite-dimensional $k$-algebras over a
field $k$. Assume that two bimodules $_{A}M _{B}$ and $_{B}N_{A}$
define a stable equivalence of Morita type between $A$ and $B$. Let
$X$ be an $A$-module which is a generator for $A$-module category.
Then, for any finite admissible set $\Phi$ of natural numbers, there
is a stable equivalence of Morita type between $\E^{\Phi}_{A}(X)$
and $\E^{\Phi}_{B}(N\otimes_{A}X)$. \label{thm1}
\end{Theo}

Note that if $\Phi=\{0\}$, then the above result was known in
\cite{LX3}. Thus Theorem \ref{thm1} generalizes the main result in
\cite{LX3}, and provides much more possibilities for constructing
stable equivalences of Morita type through the choices of different
$\Phi$ . Also, our proof of Theorem \ref{thm1} is different from
that in \cite{LX3}.

Next, we shall exploit certain kinds of restrictions to construct
stable equivalences of Morita type. Our result along this line is
the following theorem.

\begin{Theo}{\rm (The Restriction Procedure)}

Suppose that $A$ and $B$ are finite-dimensional $k$-algebras over a
field $k$ such that neither $A$ nor $B$  has semisimple direct
summands. Further, suppose that $_AM_B$ and $_BN_A$ are bimodules
without projective bimodules as direct summands, and define a stable
equivalence of Morita type between $A$ and $B$. If $e^2=e\in A$ such
that $M\otimes_BNe\in \add(Ae)$, and if $f^2=f\in B$ such that
$\add(Bf)=\add(Ne)$, then the bimodules $eMf$ and $fNe$ define a
stable equivalence of Morita type between $eAe$ and $fBf$. Moreover,
if we define $\Lambda=\End_{eAe}(eA)$, $\Gamma=\End_{fBf}(fB)$,
$N'=\Hom_{fBf}((fB)_{\Gamma},fNe\otimes_{eAe}(eA)_\Lambda)$ and
$M'=\Hom_{eAe}((eA)_\Lambda,eMf\otimes_{fBf}(fB)_{\Gamma})$, then
$_\Gamma N'_\Lambda$ and $_\Lambda M'_\Gamma$ define a stable
equivalence of Morita type between $\Lambda$ and $\Gamma$.
\label{thm2}
\end{Theo}

In fact, under the assumptions of Theorem \ref{thm2}, we may have a
more general formulation, namely, for any finite admissible set
$\Phi$ of natural numbers and for any $eAe$-module $X$, the
$\Phi$-Auslander-Yoneda algebras of $eAe\oplus X$ and $fBf\oplus
fNe\otimes_{eAe}X$ are stably equivalent of Morita type. This is a
consequence of Theorem \ref{thm1} and Theorem \ref{thm2}.

Also, from Theorem \ref{thm1} and Theorem \ref{thm2} we have the
following characterization of stable equivalences of Morita type for
representation-finite algebras as well as for self-injective
algebras.

\begin {Koro}\label{koro1} Suppose that $A$ and $B$ are finite-dimensional
$k$-algebras over a field $k$ such that neither $A$ nor $B$ has
semisimple direct summands.

$(1)$ Assume further that $A$ and $B$ are self-injective. Let $X$ be
an $A$-module and let $Y$ be a $B$-module. If there is a finite
admissible set $\Phi$ of natural numbers such that
$\E^{\Phi}_A(A\oplus X)$ and $\E^{\Phi}_B(B\oplus Y)$ are stably
equivalent of Morita type, then, for \emph{any} finite admissible
set $\Psi$ of natural numbers, the algebras $\E^{\Psi}_A(A\oplus X)$
and $\E^{\Psi}_B(B\oplus Y)$ are stably equivalent of Morita type.

$(2)$ Assume additionally that $A$ and $B$ are
representation-finite. Then $A$ and $B$ are stably equivalent of
Morita type if and only if so are their Auslander algebras.
\end{Koro}

Note that the ``only if " part of Corollary \ref{koro1} (2) follows
from \cite{LX3}.

Of course, there are many important classes of algebras which are of
the form $\End_A(A\oplus Y)$ with $A$ self-injective and $Y$ an
$A$-module. For example, Schur algebras or $q$-Schur algebras. Thus,
as a consequence of Corollary \ref{koro1}, we know that the global
dimension of $\End_{k[S_n]}(k[S_n]\oplus \Omega^i(Y))$ is finite for
$i\in \mathbb Z$, where $k[S_n]$ is the group algebra of the
symmetric group $S_n$, $Y$ the direct sum of non-projective
indecomposable Young modules, and $\Omega$ the usual syzygy
operator.

As another byproduct of our considerations in this paper, we can
construct a family of derived equivalent algebras with certain
special properties.

\begin{Koro} Suppose that $k$ is a field with a non-zero element that is not a root of unity.
Then, there is an infinite series of $k$-algebras of the same
dimension such that they have the same dominant and global
dimensions, and are all derived equivalent, but pairwise not stably
equivalent of Morita type. \label{koro2}
\end{Koro}

The contents of this paper are organized as follows. In Section
\ref{pre}, we fix notations and prepare some basic facts for our
proofs. In Section \ref{induction} and Section \ref{restriction}, we
prove our main results, Theorem \ref{thm1} and Theorem \ref{thm2},
as well as Corollary \ref{koro1}(2), respectively. In Section
\ref{selfinjective}, we concentrate our consideration on
self-injective algebras, and establish some applications of our main
results. In particular, in this section we prove Corollary
\ref{koro1}(1) and supply a sufficient condition, which is used in
Section \ref{family}, to verify when two algebras are not stably
equivalent of Morita type. In Section \ref{family}, we apply our
results in the previous sections to Liu-Schulz algebras and give a
proof of Corollary \ref{koro2} which answers a question by Thorsten
Holm.

C.C.Xi, the corresponding author, is partially supported by the
Fundamental Research Funds for the Central Universities (2009SD-17),
while H.X.Chen is supported by the Doctor Funds of the Beijing
Normal University. This revision of the first draft was partially
done when C.C.Xi visited the Chern Institute of Mathematics,
Tianjin, China, in July, 2010, he would like to thank Professor
Chengming Bai at the Nankai University for his warm hospitality. The
authors thank Yuming Liu for some helpful discussions on the
subject.

\section{Preliminaries\label{pre}}
In this section, we shall fix some notations, and recall some
definitions and basic results which are needed in the proofs of our
main results.

\subsection{Some conventions and homological facts}
Throughout this paper, $k$ stands for a fixed field. All categories
and functors will be $k$-categories and $k$-functors, respectively.
Unless stated otherwise, all algebras considered are
finite-dimensional $k$-algebras, and all modules are finitely
generated left modules.

Let $\mathcal C$ be a category. Given two morphisms $f: X\to Y$ and
$g: Y\to Z$ in $\mathcal C$, we denote the composition of $f$ and
$g$ by $fg$ which is a morphism from $X$ to $Z$, while we denote the
composition of a functor $F:\mathcal {C}\to \mathcal{D}$ between
categories $\mathcal C$ and $\mathcal D$ with a functor $G:
\mathcal{D}\ra \mathcal{E}$ between categories $\mathcal D$ and
$\mathcal E$ by $GF$ which is a functor from $\mathcal C$ to
$\mathcal E$.

If $\cal C$ is an additive category and $X$ is an object in
$\mathcal{C}$, we denote by $\add(X)$ the full subcategory of
$\mathcal{C}$ consisting of all direct summands of  direct sums of
finitely many copies of $X$. The object $X$ is called an additive
generator for $\mathcal C$ if add($X$) = $\mathcal C$.

Let $A$ be an algebra. We denote by $A\modcat$ the category of all
$A$-modules, by $A\pmodcat$ (respectively, $A\imodcat$) the full
subcategory of $A\modcat$ consisting of projective (respectively,
injective) modules,  by $D$ the usual $k$-duality $\Hom_{k}(-, k)$,
and by $\nu_{A}$ the Nakayama functor $D\Hom_{A}(-,\,_{A}A)$ of $A$.
Note that $\nu_A$ is an equivalence from $A\pmodcat$ to $A\imodcat$
with the inverse $\Hom_A(D(A),-)$. We denote the global and dominant
dimensions of $A$ by gl.dim$(A)$ and dom.dim$(A)$, respectively.

As usual, by $\Db{A}$ we denote the bounded derived category of
complexes over $A$-mod. It is known that $A$-mod is fully embedded
in $\Db{A}$ and that $\Hom_{\Db{A}}(X,Y[i])\simeq \Ext^i_A(X,Y)$ for
all $i\ge 0$ and all $A$-modules $X$ and $Y$.

Let $X$ be an $A$-module. We denote by $\Omega^i_A(X)$ the $i$-th
syzygy, by $\soc(X)$ the socle, and by rad$(X)$ the Jacobson radical
of $X$.

Let $X$ be an additive generator for $A$-mod. The endomorphism
algebra of $X$ is called the Auslander algebra of $A$. This algebra
is, up to Morita equivalence, uniquely determined by $A$. Note that
Auslander algebras can be described by two homological properties:
An algebra $A$ is an Auslander algebra if gl.dim($A$)$\leq 2\le
\dm(A)$.

An $A$-module $X$ is called a generator for $A$-mod if
$\add(_AA)\subseteq \add(X)$; a cogenerator for $A$-mod if
$\add(D(A_A))\subseteq \add(X)$, and a generator-cogenerator if it
is both a generator and a cogenerator for $A$-mod. Clearly, an
additive generator for $A$-mod is a generator-cogenerator for
$A$-mod. But the converse is not true in general.

Let $T$ be an arbitrary $A$-module, and let $B$ be the endomorphism
algebra of $T$. We consider the following full subcategories of
$A$-mod related to $T$.
$$\begin{array}{rl} Gen(_{A}T) :=&\{X\in A\modcat \mid \mbox{ there is
a surjective homomorphism from  } T^{m} \mbox{ to } X \mbox{with } m\ge 1 \}.\\
Pre(_{A}T) :=&\{X\in A\modcat \mid \mbox{there is an exact sequence
}T_{1}\rightarrow T_{0}\rightarrow X \mbox{ with all } T_i\in
\add(_{A}T) \}.\\
App(_AT) :=&\{X\in A\modcat \mid  \mbox{there is a homomorphism } g:
T_0\ra X \mbox{ with  } T_0\in \add(_AT) \mbox{ such that } \\
 & \qquad \qquad \qquad \mbox{Ker}(g)\in Gen(_{A}T) \mbox{ and  } \Hom_A(T',g) \mbox{ is surjective for  } T'\in
\add(T)\}.
 \end{array}$$

The following lemma is known, for a proof, we refer, for example, to
\cite[Lemma 2.1]{C5}.

\begin{Lem} Let $T$ be an $A$-module and $B=\End(_AT)$. Let $X$ be an arbitrary $A$-module.
Then:

$(1)$ If $Y$ is a right $B$-module, then the natural homomorphism
$\delta$: $Y\otimes_B \Hom_{A}(T,X)\ra
 \Hom_{B}(\Hom_{A}(X,T), Y)$, given by $y\otimes f \mapsto \delta_{y\otimes f}$
with $\delta_{y\otimes f}(g)= y(fg)$ for $y\in Y, f\in \Hom_A(T,X),
g\in \Hom_A(X,T)$, is an isomorphism if $X\in
 \emph{add}(_{A}T)$.

$(2)$ If $X'\in \add (_AT)$, or $ X \in \add(_AT)$, then the
composition map $\mu: \Hom_A(X',T)\otimes_B\Hom_A(T,X)\rightarrow
\Hom_A(X',X)$ given by $f\otimes_B g\mapsto fg $ is bijective.

$(3)$ If $X\in Gen(_{A}T)$, then the evaluation map $e_{X}:
T\otimes_{B}\Hom_A(T,X)\rightarrow X$ is surjective. If $X\in
App(_{A}T)$, then $e_{X}$ is bijective. Conversely, if $e_{X}$ is
bijective, then $X\in App(_{A}T)$. \label{2.3}
\end{Lem}

The next lemma is taken from \cite[Lemma 2.1]{C3}, which can also be
verified directly.

\begin{Lem}$\cite{C3}$ $(1)$ Let $A, B, C$
and $E$ be $k$-algebras, and let $_{A}X_{B}$ and $_{B}Y_{E}$ be
bimodules with $X_{B}$ projective. Put $X^{*}=\Hom_{B}(X,B)$. Then
the natural homomorphism $\varphi: {}_AX\otimes_{B}Y_{E}\rightarrow
\Hom_{B}(_BX^*_A, {}_BY_E)$, defined by $x\otimes y\mapsto
\varphi_{x\otimes y}$, where $\varphi_{x\otimes y}(f)=(xf)y$ for
$x\in X, y\in Y$ and $f\in X^{*}$, is an isomorphism of
$A$-$E$-bimodules, where the image of $x$ under $f$ is denoted by
$xf$.

$(2)$ In the situation $(_EP_A, {}_CX_{B}, {}_AU_B)$, if $P_A$ is
projective, or if $X_B$ is projective, then
$_EP\otimes_A\Hom_B(_CX_B, {}_AU_B)\simeq \Hom_B(_CX_B,
{}_EP\otimes_{A}U_B)$ as $E$-$C$-bimodules. Dually, in the situation
$(_AP_E, {}_BX_C, {}_BU_A)$, if $ _AP$ is projective, or if $_{B}X$
is projective, then $\Hom_B(_BX_C, {}_BU_A)\otimes_AP_E\simeq
\Hom_B(_BX_C, {}_BU\otimes_AP_E)$ as $C$-$E$-bimodules.\label{2.4}
\end{Lem}

The following is a well-known result due to Auslander (for example,
see \cite[Proposition 5.6, p.214]{ARS}).

\begin{Lem} Let $\Lambda$ be an Artin algebra such that $\emph{gl.dim}(\Lambda)\le 2 \le
 \emph{dom.dim}(\Lambda)$. Let $U$ be a $\Lambda$-module such that
$\add(U)$ is the full subcategory of $\Lambda\emph{-mod}$ consisting
of all projective-injective $\Lambda$-modules. Then

$(1)$ $A:= \End_{\Lambda}(U)$ is representation-finite.

$(2)$ $\Lambda$ is Morita equivalent to $\End_A(X)^{op}$, where $X$
is an additive generator for $A\emph{-mod}$. \label{2.5}
\end{Lem}

Finally, we recall the definition of $\mathcal D$-split sequences
from \cite{hx2}. For our purpose, we just restrict our attention to
module categories.

Let $\mathcal D$ be a full subcategory of $A$-mod. A short exact
sequence
$$ 0\lra X\lraf{f} M\lraf{g}Y\lra 0$$
in $A$-mod is called a $\mathcal D$-split sequence if $M\in
{\mathcal D}$, $\Hom_A(D',g)$ and $\Hom_A(f,D')$ are surjective for
every object $D'\in {\mathcal D}$.

Note that $\cal D$-split sequences were used in \cite{hx2} to
construct tilting modules of projective dimension at most one.

\subsection{Admissible sets and perforated orbit categories\label{adm}}

In \cite{HAU}, a class of algebras, called $\Phi$-Auslander-Yoneda
algebras, were introduced, which include, for example, Auslander
algebras, generalized Yoneda algebras and certain trivial
extensions.

Let $\mathbb N$ be the set of natural numbers $\{0,1,2, \cdots\}$.
Recall that a subset $\Phi$ of $\mathbb N$ is said to be admissible
provided that $0\in \Phi$ and that for any $p,q,r\in \Phi$ with
$p+q+r\in \Phi$ we have $p+q\in \Phi$ if and only if $q+r\in \Phi$.

As shown in \cite{HAU}, there are a lot of admissible subsets of
$\mathbb N$. For example, given any subset $S$ of $\mathbb N$
containing $0$, the set $\{x^m\mid x\in S\}$ is admissible for all
$m\ge 3.$

Let $\Phi$ be an admissible subset of $\mathbb N$.

Let $\mathcal C$ be a $k$-category, and let $F$ be an additive
functor from $\mathcal C$ to itself. The $(F,\Phi)$-orbit category
${\mathcal C}^{F,\Phi}$ of $\mathcal C$ is a category in which the
objects are the same as that of $\mathcal C$, and the morphism set
between two objects $X$ and $Y$ is defined to be

$$ \Hom_{{\mathcal C}^{F,\Phi}}(X,Y):=\displaystyle \bigoplus_{i\in \Phi} \Hom_{\mathcal C}(X, F^iY)\in k\Modcat,$$
and the composition is defined in an obvious way. Since $\Phi$ is
admissible, ${\mathcal C}^{F,\Phi}$ is an additive $k$-category. In
particular, $\Hom_{{\mathcal C}^{F,\Phi}}(X, X)$ is a $k$-algebra
(which may not be finite-dimensional), and $\Hom_{{\mathcal
C}^{F,\Phi}_A}(X, Y)$ is an $\Hom_{{\mathcal C}^{F,\Phi}}(X,
X)$-$\Hom_{{\mathcal C}^{F,\Phi}}(Y, Y)$-bimodule. For more details,
we refer the reader to \cite{HAU}. In this paper, the category
${\mathcal C}^{F,\Phi}$ is simply called a perforated orbit
category, and the algebra $\Hom_{{\mathcal C}^{F,\Phi}}(X, X)$ is
called the perforated Yoneda algebra of $X$ without mentioning $F$
and $\Phi$.

In case $\mathcal C$ is the bounded derived category $\Db{A}$ with
$A$ a $k$-algebra, and $F$ is the shift functor [1] of $\Db{A}$, we
denote simply by ${\mathcal E}^{\Phi}_A$ the $(F,\Phi)$-orbit
category ${\mathcal C}^{F,\Phi}$, by $\E_{A}^{\Phi}(X, Y) $ the set
$\Hom_{{\mathcal E}^{\Phi}_A}(X,Y)$, and by $\E_A^{\Phi}(X)$ the
endomorphism algebra $\Hom_{{\mathcal E}^{\Phi}_A}(X, X)$ of $X$ in
${\mathcal E}^{\Phi}_A$. It is called $\Phi$-Auslander-Yoneda
algebra of $X$. Note that each element in $\E_A^{\Phi}(X,Y)$ can be
written as $(f_i)_{i\in\Phi}$ with $f_j\in \Hom_{\Db{A}}(X, Y[j])$.
The composition of morphisms in${\mathcal E}^{\Phi}_A$ can be
interpreted as follows: for each triple $(X, Y, Z)$ of objects in
$\Db{A}$,

$$\E_{A}^{\Phi}(X, Y) \times \E_A^{\Phi}(Y, Z) \lra \E_{A}^{\Phi}(X,Z)$$
$$\big((f_u)_{u\in \Phi}, (g_v)_{v\in\Phi}\big) \mapsto (h_i)_{i\in \Phi},$$
where
$$ h_i:=\sum _{{u, v\in \Phi}\atop{u+v=i}}\,{f_u(g_v[u])}$$
for each $i\in\Phi$. Clearly, if $\Phi$ is finite, then
$\E_A^{\Phi}(X,Y)$ is finite-dimensional for all $X,Y\in A\modcat$.

Now, let us state some elementary properties of the Hom-functor
$\E^{\Phi}_A(X,-)$.

\begin{Lem}  Suppose that $A$ is an algebra,  that $X$ is
an $A$-module, and that $\Phi$ is a finite admissible subset of
$\mathbb{N}$.

$(1)$ Let $\add_A^{\Phi}(X)$ stand for  the full subcategory of
${\mathcal E}_A^{\Phi}$ consisting of objects in $\add(_AX)$. Then
the Hom-functor \,$\E_{A}^{\Phi}(X, -) : \add_A^{\Phi}(X)\lra
\E_{A}^{\Phi}(X)\pmodcat $ \,is an equivalence of categories;

$(2)$ Let $B$ be a $k$-algebra, and let $P$ be a $B$-$A$-bimodule
such that $P_A$ is projective. Then there is a canonical algebra
homomorphism $\alpha_P: \E^{\Phi}_A(X)\lra \E^{\Phi}_B(P\otimes_AX)$
defined by $(f_i)_{i\in \Phi}\mapsto (P\otimes_Af_i)_{i\in \Phi}$
for $(f_i)_{i\in \Phi}\in \E^{\Phi}_A(X)$. Thus every left (or
right) $\E^{\Phi}_B(P\otimes_AX)$-module can be regarded as a left
(or right) $\E^{\Phi}_A(X)$-module via $\alpha_P$. \label{lem1}
\end{Lem}

The following homological result plays an important role in proving
 Theorem \ref{thm1}.

\begin{Lem}
Suppose that  $A, B$ and $C$  are $k$-algebras. Let $_AX$ be a
module, and let $_AY_B$ and ${}_BP_C$ be bimodules with $_BP$
projective. Then, for each $i\geq 0$, we have $\Ext^{i}_{A}(X,
Y\otimes_B P_C)\simeq \Ext^i_A(X, Y)\otimes_B P_C $\, as
$C^{op}$-modules. Moreover, for each  admissible subset $\Phi$ of
$\mathbb{N}$, we have $\E^{\Phi}_{A}(X, Y\otimes_B P_C)\simeq
\E^{\Phi}_A(X, Y)\otimes_B P_C $ as
$\E^{\Phi}_{A}(X)\mbox{-}C$-bimodules. \label{lem3}
\end{Lem}

{\bf Proof}. First, let us recall the Yoneda product. Assume that
$U, V \mbox{and}\; W $ are $A$-modules. Fix a minimal projective
resolution $\cpx{P}_U$ of $_AU$: \;$$\cdots \lra P^n\lraf{d^n}
P^{n-1}\lra\cdots \lra P^1\lraf{d^1} P^0 \lraf{d^0} U\lra 0,$$ with
all $P^i$ projective. If $g: U\ra V$ is a homomorphism, then there
is a lifting of $g$, which is a chain map $\cpx{g}:\cpx{P}_U\ra
\cpx{P}_V$. Then, for each $i\geq 1$, we have a short exact sequence
$0\to \Omega^i_A(U)\lraf{\lambda{_i}} P^{i-1}\lraf{\mu_{i}}
\Omega^{i-1}_A(U)\ra 0$, which gives rise to a right exact sequence
of $k$-modules
$$\Hom_{A}(P^{i-1}, V) \lraf{(\lambda{_i})_{\ast}}
\Hom{_A}(\Omega^{i}_{A}(U), V)\lra \Ext^{i}_{A}(U, V)\lra 0.$$ Hence
each element of $\Ext^{i}_{A}(U, V)$ can be regarded as a
homomorphism in $\Hom{_A}(\Omega^{i}_{A}(U), V)$ modulo the subspace
of $\Hom{_A}(\Omega^{i}_{A}(U), V)$ generated by all homomorphisms
that factorize through $\lambda_{i}$, where $i\geq 0$ and
$P^{-1}:=0$. In what follows, we denote the image of $f\in
\Hom{_A}(\Omega^{i}_{A}(U), V)$ by $\ol{f}\in \Ext^{i}_{A}(U, V)$.

Given  $i, j \in \mathbb{N}, f_i\in \Hom{_A}(\Omega^{i}_{A}(U), V)$
and $g_j\in \Hom{_A}(\Omega^{j}_{A}(V), W)$, we know that the Yoneda
product $\mu: \Ext^{i}_{A}(U, V)\otimes_k\Ext^{j}_{A}(V, W)\to
\Ext^{i+j}_{A}(U, W)$ can be presented by
$\ol{f_i}\otimes_k\ol{g_j}\mapsto \ol{\Omega_A^j(f_i)g_j}$, where
$\Omega_A^j(f_i)$ is the $j$-th term of a lifting of $f_i$. Note
that the Yoneda product is independent of the choice of a lifting of
$f_i$.

By Lemma \ref{2.4}(2), for each
 $_AW$, there is a natural $C^{op}$-module isomorphism $\theta _W: \Hom{_A}(W, Y)\otimes
_B P_C\to \Hom{_A}(W, Y\otimes_BP_C)$ defined by $\theta_W(f\otimes
p)(w)= f(w)\otimes p$ for\,  $f\in \Hom{_A}(W, Y), p\in P,
\;\mbox{and} \; w\in W$. In other words, we have a natural
equivalence $\theta: \Hom{_A}(-, Y)\otimes _B P_C\simeq \Hom{_A}(-,
Y\otimes_BP_C)$ of functors from $ A\modcat$ to $C^{op}\modcat$. Let
$$\cdots\lra Q^i\lra Q^{i-1}\lra\cdots\lra Q^1\lra Q^0\lra X\lra 0$$ be a
minimal projective resolution of $_A X$. Then, by definition, we
have a right exact sequence of $k$-modules $$\Hom_{A}(Q^{i-1},
Y)\lra \Hom{_A}(\Omega^{i}_{A}(X), Y)\lra\Ext^{i}_{A}(X, Y)\lra 0.$$
Since $_BP$ is projective, the following diagram is exact and
commutative for $i\geq 0$:
$$\xymatrix{
\Hom{_A}(Q^{i-1}, Y)\otimes_B P_C\ar@{}[dr] | {\displaystyle }
\ar[d]^-{\wr}_-{\theta_{Q^{i-1}}}\ar[r]&\Hom{_A}(\Omega_A^i(X),
Y)\otimes_BP_C\ar[d]^-{\wr}_-{\theta_{\Omega_A^i(X)}}\ar[r]
&\Ext^{i}_{A}(X, Y)\otimes_BP_C\ar[r]\ar@{-->}_{\varphi_i}[d]& 0\\
\Hom{_A}(Q^{i-1}, Y\otimes_B P_C)\ar[r]&\Hom{_A}(\Omega_A^i(X),
Y\otimes_BP_C) \ar[r] &\Ext^{i}_{A}(X, Y\otimes_BP_C)\ar[r]& 0 ,}
$$ where we set $Q^{-1}:=0$.
This induces an isomorphism $\varphi_i: \Ext^{i}_{A}(X,
Y)\otimes_BP_C\to \Ext^{i}_{A}(X, Y\otimes_BP_C)$ defined by
$\ol{f_i}\otimes p\mapsto \ol{\theta_{\Omega_A^i(X)}(f_i\otimes
p)}$, where $f_i\in \Hom{_A}(\Omega_A^i(X), Y)\text{\,and\,} p\in
P$. Clearly, $\varphi_i$ is a  $C^{op}$-homomorphism for each $i\geq
0$. Thus the first part of Lemma \ref{lem3} is proved.

Second, for each admissible subset $\Phi$ of $\mathbb{N}$, we define
a map $\varphi_{\Phi}: \E^{\Phi}_A(X, Y)\otimes_B P_C \to
\E^{\Phi}_{A}(X, Y\otimes_B P_C)$ by $(\ol{f_i})\otimes p\mapsto
(\varphi_i(\ol{f_i}\otimes p))$, where $p\in P$, and $f_i\in
\Hom{_A}(\Omega_A^i(X), Y) $ with $i\in \Phi$. By the above
discussion, we know that  $\varphi_{\Phi}$ is an isomorphism of
$C^{op}$-modules. In order to prove that $\varphi_{\Phi}$ is an
isomorphism of $\E^{\Phi}_{A}(X)\mbox{-}C$-bimodules, it suffices to
show that $\varphi_{\Phi}$ is an isomorphism of left
$\E_A^{\Phi}(X)$-modules, or equivalently, we have to check that the
following diagram commutes for $i, j\in \Phi $ with $ i+j\in \Phi$:

$$
\xymatrix{ \Ext^{i}_{A}(X, X)\otimes_{k}\Ext^{j}_{A}(X,
Y)\otimes_BP_C \ar[d]^-{\mu\otimes 1} \ar[r]^-{1\otimes \varphi_j }
& \Ext^{i}_{A}(X,
X)\otimes_{k}\Ext^{j}_{A}(X, Y\otimes_BP_C)\ar[d]^-{\mu}\\
\Ext^{i+j}_{A}(X, Y)\otimes_B P_C\ar[r]^-{\varphi_{i+j}}
&\Ext^{i+j}_{A}(X, Y\otimes_BP_C), }
$$
where $\mu$ is the usual Yoneda product. Let $u\in\Hom
{_A}(\Omega_A^i(X), X)$, $v\in \Hom{_A}(\Omega_A^j(X), Y)$ and $p\in
P$. Then
$$
((1\otimes \varphi_j)\, \mu)\,(u\otimes v\otimes p)=
\ol{\Omega_A^j(u)\,\theta_{\Omega_A^j(X)}\,(v\otimes p)}
\mbox{\;and\;} ((\mu\otimes 1)\,\varphi_{i+j})\,(u\otimes v\otimes
p)= \ol{\theta_{\Omega_A^{i+j}(X)}\,((\Omega_A^j(u)v)\otimes p)}.
$$
By definition, for each $x\in \Omega_A^{i+j}(X)$,  we get
$$
(\Omega_A^j(u)\,\theta_{\Omega_A^j(X)}\,(v\otimes p))(x)=
(\Omega_A^j(u)\,v)\,(x)\otimes p =
(\,\theta_{\Omega_A^{i+j}(X)}\,((\Omega_A^j(u)v)\otimes p)) (x).
$$
It follows that  $(1\otimes \varphi_j)\, \mu = (\mu\otimes
1)\,\varphi_{i+j}$. This implies that $\varphi_{\Phi}$ is an
isomorphism of $\E^{\Phi}_{A}(X)$-$C$-bimodules. Thus the proof is
completed. $\square$

\section{Inductions for stable equivalences of Morita type \label{induction}}

In this section, we shall prove Theorem \ref{thm1}. First, we recall
the definition of stable equivalences of Morita type in \cite{B}.

\begin{Def}
Let $A$ and $B$ be (arbitrary) $k$-algebras. We say that $A$ and $B$
are stably equivalent of Morita type if there is an $A$-$B$-bimodule
$_A M_B$ and a $B$-$A$-bimodule $_B N_ A$ such that

$(1)$ $M$ and $N$ are projective as one-sided modules, and

$(2)$  $M\otimes_{B}N\simeq {A\oplus P}$ as $A$-$A$-bimodules for
some projective $A$-$A$-bimodule $P$, and $N\otimes_{A}M \simeq
{B\oplus Q}$ as $B$-$B$-bimodules for some projective
$B$-$B$-bimodule $Q$. \label{stm}
\end{Def}

In this case, we say that $M$ and $N$ define a stable equivalence of
Morita type between  $A$ and $B$. Moreover, we have two exact
functors $T_N:=N\otimes_A- : A\modcat \to B\modcat$ and $T_M
:=M\otimes_B-: B\modcat \to A\modcat$. Similarly, the bimodules $P$
and $Q$ define two exact functors $T_P$ and $T_Q$, respectively.
Note that the images of $T_P$ and $T_Q$ consist of projective
modules.

From now on, we assume that $A, B, M, N, P$ and $Q$ are fixed as in
Definition \ref{stm}, and that $X$ is a generator for $A\modcat$.
Moreover, we fix a finite admissible subset $\Phi$ of $\mathbb{N}$,
and define $\Lambda:=\E^{\Phi}_{A}(X)$ and
$\Gamma:=\E^{\Phi}_{B}(N\otimes_A X)$.

Since the functors $T_N$ and $T_M$ are exact, they preserve
acyclicity, and can be extended to  triangle functors $T_N': \Db
A\to\Db B$ and $T_M': \Db B\to\Db A$, respectively. Furthermore,
$T'_N$ and $T'_M$ induce canonically two functors $F: {\mathcal
E}_A^{\Phi}\to {\mathcal E}_B^{\Phi}$ and $G: {\mathcal
E}_B^{\Phi}\to {\mathcal E}_A^{\Phi}$, respectively. More precisely,
if $\cpx{X}\in \Db{A}$, then $F(\cpx{X}):= ( N\otimes_AX^i,
N\otimes_Ad_X^i)$, and if $f:=(f_j)_{j\in \Phi}\in \Hom_{{\mathcal
E}^{\Phi}_A}(\cpx{X},\cpx{Y})$ with $\cpx{Y}\in \Db{A}$, then
$F(f):=(N\otimes_Af_j)_{j\in \Phi}\in \Hom_{{\mathcal
E}^{\Phi}_B}(F(\cpx{X}),F(\cpx{Y}))$. Similarly, we define the
functor $G$.

The functor $F$ gives rise to a canonical algebra homomorphism
$\alpha_N: \E_A^{\Phi}(\cpx{X})\to\E_B^{\Phi}(F(\cpx{X}))$ for each
object $\cpx{X}\in \Db{A}$. In particular, for any
$\cpx{Z}\in\Db{B}$, we can regard $\E_B^{\Phi}(\cpx{Z}, F(\cpx{X}))$
as  an $\E_B^{\Phi}(\cpx{Z})$-$\E_A^{\Phi}(\cpx{X})$-bimodule via
$\alpha_N$. Note that  the homomorphism $\alpha_N$ coincides with
the one defined in Lemma \ref{lem1}, when $\cpx{X}$ is an
$A$-module.

\medskip
{\bf Proof of Theorem \ref{thm1}.}  We define $U:=\E^{\Phi}_{A}(X,
T_M (N\otimes_A X))$ and $V:=\E^{\Phi}_{B}(N\otimes_A X, T_N (X))$.
In the following we shall prove that $U$ and $V$ define a stable
equivalence of Morita type between $\Lambda$ and $\Gamma$.

First, we endow $U$ with a right $\Gamma$-module structure by
$u\cdot \gamma:=uG(\gamma)$ for $u\in U$ and $\gamma\in \Gamma$, and
endow $V$ with a right $\Lambda$-module structure by
$v\cdot\lambda:=vF(\lambda)$ for $v\in V$ and $\lambda\in \Lambda$.
Then, $U$ becomes  a $\Lambda$-$\Gamma$-bimodule, and  $V$ becomes a
$\Gamma$-$\Lambda$-bimodule.

By definition, we know $V=\Gamma$, and it is a projective left
$\Gamma$-module. Note that $_AX$ is a generator and  the images of
$T_P$ consists of projective modules. We conclude  that $T_M
(N\otimes_A X) = M\otimes_B (N\otimes_A X)\simeq X\oplus P\otimes_A
X\in\add(X)$. Thus $U$ is projective as a left $\Lambda$-module by
Lemma \ref{lem1}.

(1) $U \otimes_\Gamma V$, as a $\Lambda$-$\Lambda$-bimodule,
satisfies the condition (2) in Definition \ref{stm}.

Indeed, we write $ W:=\E^{\Phi}_{A}(X, (T_MT_N)(X))$, and define a
right $\Lambda$-module structure on $W$ by $w\cdot
\lambda':=w(GF)(\lambda')$ for  $w\in W$ and $\lambda'\in \Lambda$.
Then $W$ becomes a $\Lambda$-$\Lambda$-bimodule. Note  that  there
is a natural $\Lambda$-module isomorphism $\varphi: U \otimes_\Gamma
V \ra W$ defined by $x\otimes y\mapsto xG(y)$ for $x\in U$ and $y\in
V$. We claim that $\varphi$ is an isomorphism of
$\Lambda$-$\Lambda$-bimodules. In fact, it suffices to show that
$\varphi$ respects the structure of right $\Lambda$-modules.
However, this follows immediately from a verification: for $c\in U,
d\in V$ and $a\in \Lambda$, we have
$$
\varphi((c\otimes d)\cdot a)=\varphi(c\otimes (d F(a)))=cG(d
F(a))=cG(d)(GF)(a) =\varphi(c\otimes d)\cdot a.
$$
Combining this bimodule isomorphism $\varphi$ with Lemma \ref{lem1},
we get the following isomorphisms of $\Lambda$-$\Lambda$-bimodules:
 $$ (*)\qquad U \otimes_\Gamma V \simeq
\E^{\Phi}_{A}(X, (T_MT_N)(X))\simeq \E^{\Phi}_{A}(X, X)\oplus
\E^{\Phi}_{A}(X, P\otimes_A X)=\Lambda\oplus \E^{\Phi}_{A}(X,
P\otimes_A X),
$$
where  the second isomorphism follows  from  $M\otimes_{B}N\simeq
{A\oplus P}$ as $A$-$A$-bimodules, and where the right
$\Lambda$-module structure on $\E^{\Phi}_{A}(X, P\otimes_A X)$ is
induced by the canonical algebra homomorphism $\Lambda \to
\E^{\Phi}_{A}( P\otimes_A X)$, which sends $(f_i)_{i\in \Phi}$ in
$\Lambda$ to $(P\otimes_Af_i)_{i\in \Phi}$ (see Lemma \ref{lem1}
(2)).

Now, we show that $\E^{\Phi}_{A}(X, P\otimes_A\,X)$ is a projective
$\Lambda$-$ \Lambda$-bimodule. In fact, since $P\in\add(_AA\otimes_k
A_A)$, we conclude that $\E^{\Phi}_{A}(X,
P\otimes_A\,X)\in\add(\E^{\Phi}_{A}(X, (A\otimes_k A)
\otimes_A\,X))$. Thus, it is sufficient to prove that
 $\E^{\Phi}_{A}(X, (A\otimes_k A) \otimes_A\,X)$ is a
projective $\Lambda$-$ \Lambda$-bimodule. For this purpose, we first
note that the right $\Lambda$-module structure on $\E^{\Phi}_{A}(X,
(A\otimes_k A) \otimes_A\,X)$ is induced by the canonical algebra
homomorphism  $\alpha_{A\otimes_kA}:\Lambda\to
\E^{\Phi}_{A}((A\otimes_k A) \otimes_AX)$, which sends
$g:=(g_i)_{i\in \Phi}$ in $\Lambda$ to
$((A\otimes_kA)\otimes_Ag_i)_{i\in \Phi}$. Clearly, $_AA\otimes_k
A\otimes_A\,X\in\add(_AA)$. It follows that $\Ext_A^{j}((A\otimes_k
A) \otimes_AX, (A\otimes_k A) \otimes_AX)=0$ for any $j>0$, and
therefore $(A\otimes_kA)\otimes_Ag_i=0$ for any $0\ne i\in\Phi$.
Thus we have $\alpha_{A\otimes_kA}(g)=(A\otimes_kA)\otimes_Ag_0$. If
$\pi:\Lambda\to\End_A(X)$ is the canonical projection  and $\mu'$ is
the canonical algebra homomorphism $\End_A(X)\ra
\End_A\big((A\otimes_k A)\otimes_AX\big)$, then
$\alpha_{A\otimes_kA}=\pi\mu'$. Thus the right $\Lambda$-module
structure on $\E^{\Phi}_{A}(X, (A\otimes_k A) \otimes_A\,X)$ is
induced by $\End_A(X)$. Similarly, from the homomorphisms
$$\Lambda=\E^{\Phi}_{A}(X)\epa {\pi}\End_A(X)\epa {\mu}
\End_A(A\otimes_k X)= \E^{\Phi}_{A}(A\otimes_k X),$$ where
$\mu:\End_A(X)\to \End_A(A\otimes_kX)$ is induced by the tensor
functor $A\otimes_k-$, we see that the right $\Lambda$-module
structure on $\E^{\Phi}_{A}(X, A\otimes_k X)$ is also induced by
$\End_A(X)$. Thus $\E^{\Phi}_{A}(X, (A\otimes_k A)
\otimes_A\,X)\,\simeq\,\E^{\Phi}_{A}(X, A\otimes_k X)$ as
$\Lambda$-$\Lambda$-bimodules. Moreover,  it follows from Lemma
\ref{lem3} that $\E^{\Phi}_{A}(X, A\otimes_k X)\,\simeq\,
 \E^{\Phi}_{A}(X, A)\otimes_k X$ as
$\Lambda$-$\End_A(X)$-bimodule. Since the $A$-module $X$ can be
regarded as a right $\Lambda$-module via the homomorphism $\pi$, we
see that $X$ is actually isomorphic to $\E^{\Phi}_{A}(A, X)$ as
right $\Lambda$-modules. Thus $ \E^{\Phi}_{A}(X, A)\otimes_k
X\,\simeq\,\E^{\Phi}_{A}(X, A) \otimes_k\E^{\Phi}_{A}(A, X)$ as
$\Lambda$-$\Lambda$-bimodules. Since $_AA\in\add(X)$, we know that
$\E^{\Phi}_{A}(X, A)$ is a projective $\Lambda$-module and
$\E^{\Phi}_{A}(A, X)$ is a projective right $\Lambda$-module. Hence
$ \E^{\Phi}_{A}(X, A)\otimes_k X$ is a projective
$\Lambda$-$\Lambda$-bimodule. This implies that $\E^{\Phi}_{A}(X,
P\otimes_A\,X)$ is a projective $\Lambda$-$ \Lambda$-bimodule.

(2) $V \otimes_\Lambda U$, as a $\Gamma$-$\Gamma$-bimodule, fulfills
the condition (2) in Definition \ref{stm}.

Let $Z:=\E^{\Phi}_{B}(N\otimes_A X, T_N T_M(N\otimes_A X))$.
Similarly, we endow $Z$ with a right $\Gamma$-module structure
defined by $z\cdot b:=z (FG)(b)$ for $z\in Z$ and $b\in \Gamma$.
Then $Z$ becomes a $\Gamma$-$\Gamma$-bimodule. Observe that, for
each $A$-module $Y$, there is a homomorphism
 $\Psi_Y: V\otimes_{\Lambda}\E^{\Phi}_{A}(X,
Y)\to \E^{\Phi}_{B}(N\otimes_A X, T_N(Y))$ of $\Gamma$-modules,
which is defined by $g\otimes h \mapsto gF(h)$ for $g\in V$ and
$h\in \E^{\Phi}_{A}(X, Y)$. This homomorphism is natural in $Y$. In
other words, $\Psi:V\otimes_{\Lambda}\E^{\Phi}_{A}(X, -)\to
\E^{\Phi}_{B}(N\otimes_A X, T_N(-))$  is a natural transformation of
functors from $A\modcat$ to $\Gamma\modcat$. Clearly, $\Psi_X$ is an
isomorphism of $\Gamma$-modules. It follows from $T_M (N\otimes_A
X)\in\add(X)$ that $\Psi{_{T_M (N\otimes_A X)}}: V \otimes_\Lambda U
\to Z$ is a $\Gamma$-isomorphism. Similarly, we can check that
$\Psi{_{T_M (N\otimes_A X)}}$ preserves the structure of right
$\Gamma$-modules. Thus $\Psi{_{T_M (N\otimes_A X)}}: V
\otimes_\Lambda U \to Z$ is an isomorphism of
$\Gamma$-$\Gamma$-bimodules, and there are the following
isomorphisms of $\Gamma$-$\Gamma$-bimodules:
$$ (**) \qquad V \otimes_\Lambda U \simeq Z\simeq
\Gamma \oplus \E^{\Phi}_{B}(N\otimes_A X, Q\otimes_B(N\otimes_A X)),
$$ where the second isomorphism is deduced from
$N\otimes_{A}M\simeq {B\oplus Q}$ as $B$-$B$-bimodules. By an
argument similar to that in the proof of (1), we can show that
$\E^{\Phi}_{B}(N\otimes_A X, Q\otimes_B(N\otimes_A X))$ is a
projective $\Gamma$-$\Gamma$-bimodule.

It remains to show that $U_{\Gamma}$ and $V_{\Lambda}$ are
projective. This is equivalent to showing that the tensor functors
$T_U:= U\otimes_\Gamma-: \Gamma\modcat \to\Lambda\modcat$ and
$T_V:=V\otimes_\Lambda-: \Lambda\modcat \to\Gamma\modcat$ are exact.
Since tensor functors are always right exact, the exactness of $T_U$
is equivalent to the property that $T_U$ preserve injective
homomorphisms of modules. Now, suppose that $f:C\to D$ is an
injective homomorphism between $\Gamma$-modules $C$ and $D$. Since
$\E^{\Phi}_{B}(N\otimes_A X, Q\otimes_B(N\otimes_A X))$ is a right
projective $\Gamma$-module, we know from $(**)$ that the composition
functor $T_VT_U$ is exact. In particular, the homomorphism
$(T_VT_U)(f):(T_VT_U)(C)\to (T_VT_U)(D)$ is injective. Let $\mu:
\Ker(T_U(f))\to T_U(C)$ be the canonical inclusion. Clearly, we have
$\mu T_U(f)=0$, which shows  $T_V(\mu
T_U(f))=T_V(\mu)(T_VT_U)(f)=0$. It follows that  $T_V(\mu)=0$ and
$(T_UT_V)(\mu)=0$. By ($*$), we get $\mu=0$, which implies that the
homomorphism  $T_U(f)$ is injective. Hence $T_U$ preserves injective
homomorphisms. Similarly,  we can show that $T_V$ preserves
injective homomorphisms, too. Consequently,  $U_\Gamma$ and
$V_\Lambda$ are projective.

Thus, the bimodules $U$ and $V$ define a stable equivalence of
Morita type between $\Lambda$ and $\Gamma$. This finishes the proof
of Theorem \ref{thm1}. $\square$

\smallskip
{\it Remarks.} (1) If we take $\Phi=\{0\}$ in Theorem \ref{thm1},
then we get \cite[Theorem 1.1]{LX3}. If we assume that $A$ is a
self-injective algebra, then we get a stable equivalence of Morita
type between $\E^{\Phi}_A(A\oplus X)$ and $\E^{\Phi}_A(A\oplus
\Omega_A^i(X))$ for any $A$-module $X$, any finite admissible subset
$\Phi$ of $\mathbb N$, and any integer $i\in {\mathbb Z}$. This
follows from Theorem \ref{thm1} and the fact that $\Omega_A$
provides a stable equivalence of Morita type between $A$ and itself
if $A$ is self-injective. Thus we re-obtain the stable equivalence
of \cite[Corollary 3.14]{HAU}.

(2) Since stable equivalences of Morita type preserve the global,
dominant and finitistic dimensions of algebras, Theorem \ref{thm1}
asserts actually also that these dimensions are equal for algebras
$\E_A^{\Phi}(X)$ and $\E_B^{\Phi}(N\otimes_AX)$.

Many important classes of algebras are of the form $\End_A(A\oplus
X)$ with $A$ a self-injective algebra. From the above remarks (see
also \cite[Corollary 3.14]{HAU}), we may get a series of algebras
which are stably equivalent of Morita type to Schur algebras. For
unexplained terminology in the next corollary, we refer the reader
to \cite{Gr}.

\begin{Koro} Suppose that $k$ is an algebraically closed field. Let $S_n$ be the symmetric group of degree $n$.
We denote by $Y$ the direct sum of all non-projective Young modules
over the group algebra $k[S_n]$ of $S_n$. Then,

$(1)$ for  every finite admissible subset $\Phi$ of $\mathbb N$, the
algebras $\E^{\Phi}_{k[S_n]}(k[S_n]\oplus Y)$ and
$\E^{\Phi}_{k[S_n]}(k[S_n]\oplus \Omega^i(Y))$ are stably equivalent
of Morita type for all $i\in \mathbb Z$.

$(2)$ All algebras $\End_{k[S_n]}(k[S_n]\oplus \Omega^i(Y))$ are
stably equivalent of Morita type to the Schur algebra $S_k(n,n)$. In
particular, ${\rm gl.dim}\big(\End_{k[S_n]}(k[S_n]\oplus
\Omega^i(Y)\big)< \infty$ for all $i\in {\mathbb Z}.$
\end{Koro}

\section{Restrictions for stable equivalences of Morita type\label{restriction}}

In this section, we shall consider the general question of how to
transfer stable equivalences of Morita type between algebras $A$ and
$B$ over a field to the ones between $eAe$ and $fBf$, where $e$ and
$f$ are idempotent elements in $A$ and $B$, respectively. In
particular, we shall prove Theorem \ref{thm2} in this section.

Before we start with our proof of Theorem \ref{thm2}, we state the
following facts, which are essentially known in literature. However,
we would like to collect them together as a lemma for the
convenience of the reader.

\begin{Lem} Suppose that  $A$ and $B$
are $k$-algebras without semisimple direct summands. Assume that
$_AM_B$ and $_BN_A$  define a stable equivalence of Morita type
between $A$ and $B$, and that $M$ and $N$ do not have any projective
bimodules as direct summands. Then,

$(1)$ there are isomorphisms of bimodule: $N\simeq
\Hom_{A}(M,A)\simeq \Hom_{B}(M,B)$ and $M\simeq \Hom_{A}(N,A)\simeq
\Hom_{B}(N,B)$.

$(2)$ Both $(N \otimes_{A} -, M \otimes_{B} - )$ and $( M
\otimes_{B} -, N\otimes_{A}-)$ are adjoint pairs of functors.

$(3)$ There are isomorphisms of bimodules: $P\simeq\Hom_A(P,A)$ and
$Q\simeq\Hom_B(Q,B)$, where $P$ and $Q$ are bimodules defined in
Definition \ref{stm}. Moreover, the bimodules $_AP_A$ and $_BQ_B$
are projective-injective.

$(4)$ If $_{A}I$ is injective, then so is $N\otimes_A I$.
 \label{2.6}

\end{Lem}
{\it Proof.} $(1)$  Note that, if $M$ and $N$ are indecomposable
bimodules, then all the statements in Lemma \ref{2.6} have been
proved in \cite[Theorem 2.7, Corollary 3.1, Lemma 3.2]{DM} under the
hypothesis of separability on the semisimple quotient algebras
$A/\rad(A)$ and $B/\rad(B)$. One can check that they are still valid
without the hypothesis of separability condition. In the following,
we shall use \cite[Theorem 2.2]{L1} to show Lemma \ref{2.6} under
the weaker assumption that $M$ and $N$ do not have any projective
bimodules as direct summands.

Since $A$ and $B$ are stably equivalent of Morita type and do not
have any semisimple direct summands, it follows from
\cite[Proposition 2.1]{L1} that $A$ and $B$ have the same number of
indecomposable direct summands (of two-sided ideals). Suppose that
$A=A_1\times A_2\times \cdots \times A_n$ and $B=B_1\times B_2\times
\cdots\times B_n$, where all $A_i$ and all $B_i$ themselves are
indecomposable algebras. By the proof of \cite[Theorem 2.2]{L1}, we
know that, up to suitable reordering, for each $1\leq i\leq n$,
there is an $A_i$-${B_i}$-bimodule $M_i$ and a $B_i$-$A_i$-bimodule
$N_i$ such that $M_i$ and $N_i$ are direct summands of $M$ and $N$
as bimodules, respectively, and that $M_i$ and $N_i$ define a stable
equivalence of Morita type between $A_i$ and $B_i$. Set
$M':=\bz_{1\leq i\leq n} M_i$ and $ N':=\bz_{1\leq i\leq n} N_i$.
Clearly, $M'$ and $N'$ are direct summands of $M$ and $N$,
respectively. Further, one can check directly that $M'$ and $N'$
also define a stable equivalence of Morita type between $A$ and $B$.
Since $_AM_B$ and $_BN_A$ do not have any projective bimodules as
direct summands, it follows from \cite[Lemma 4.8]{LX2} that $M\simeq
M'$ as $A$-$B$-bimodules and $N\simeq N'$ as $B$-$A$-bimodules. Note
that $A_i$ and $B_i$ are indecomposable algebras, and $M_i$ and
$N_i$ do not have any projective bimodules as direct summands. Then,
by \cite[Lemma 2.1]{DM}, we conclude that $M_i$ and $N_i$ are
indecomposable bimodules. This implies that Lemma \ref{2.6} holds
for the algebras $A_i$ and $B_i$ together with the bimodules $M_i$
and $N_i$ for each $i$. Consequently, there are isomorphisms of
$B$-$A$-bimodules: $\Hom_A(M,A)\simeq\Hom_A(\bz_{1\leq u\leq n} M_u,
\bz_{1\leq v\leq n} A_v)\simeq\bz_{1\leq u\leq n}\Hom_A(M_u,
A_u)\simeq\bz_{1\leq u\leq n}N_u\simeq N$. Similarly, we can prove
another statement in (1).

$(2)$ Note that the pair $(N \otimes_{A} -, M \otimes_{B} - )$ is an
adjoint pair of functors if and only if $_AM_B\simeq \Hom_{B}(N,B)$
as bimodules. Thus $(2)$ is a consequence of $(1)$.

$(3)$ It follows from the proof of \cite[Lemma 4.5]{C3} that the
first part of $(3)$ holds true, and that $P$ and $Q$ are injective
as one-sided modules. Furthermore, we claim that $P$ is an injective
bimodule. In fact, it suffices to show  that, for any indecomposable
direct summand $P'$ of $P$, the bimodule $_AP'_A$ is injective.
Since $_AP\in\add(A\otimes_{k}A^{op})$, there are primitive
idempotents
 $e_1$ and $e_2$ of $A$ such that
$P'\in\add(Ae_1\otimes_ke_2A)$.  This implies that $Ae_1$ and $e_2A$
are injective modules because $P'$ is injective as a one-sided
module. Thus $P'$ is an injective bimodule, and so is $P$.
Similarly, we can prove that $Q$ is injective as a bimodule.

$(4)$  We observe that there is an isomorphism of $B$-modules:
$N\otimes_A I\simeq\Hom_A(M,I)$. Since $M_B$ is projective and $_AI$
is injective, we see that $\Hom_A(M,I)$ is an injective $B$-module,
and so is $N\otimes_A I$. This completes the proof of Lemma
\ref{2.6}. $\square$

\medskip
By Lemma \ref{2.6}, we have the following corollary, which provides
examples such that the conditions of Theorem \ref{thm2} are
satisfied. Note that the last statement in Corollary \ref{2.8} below
follows also from the proof of \cite[Lemma 4.5]{C3}.

\begin{Koro} Suppose that $A$ and $B$ are $k$-algebras.
Assume that $\{e_1,\cdots, e_n\}$ and $\{f_1, \cdots, f_m\}$ are
complete sets of pairwise orthogonal primitive idempotents in $A$
and in $B$, respectively. Let $e$ be the sum of all those $e_i$ for
which $Ae_i$ is projective-injective, and let $f$ be the sum of all
those $f_j$ for which $Bf_j$ is projective-injective. If $M$ and $N$
are indecomposable bimodules that define a stable equivalence of
Morita type between $A$ and $B$, then $Ne\simeq N\otimes_{A}Ae \in
\add(Bf), Mf\simeq M\otimes_{B}Bf \in \add(Ae)$, and $Pe\in
\add(Ae)$. \label{2.8}
\end{Koro}

{\bf Proof of Theorem \ref{thm2}}. Let us remark that if $A$ and $B$
have no separable direct summands, then we may assume that $M$ and
$N$ have no non-zero projective bimodules as direct summands. In
fact, If $M=M'\oplus M''$ and $N=N'\oplus N''$ such that $M'$ and
$N'$ have no non-zero projective bimodules as direct summands, and
that $M''$ and $N''$ are projective bimodules, then it follows from
\cite[Lemma 4.8]{LX2} that $M'$ and $N'$ also define a stable
equivalence of Morita type between $A$ and $B$.

Suppose that $_AM_B$ and $_BN_A$ do not have any non-zero projective
bimodules as direct summands, and define a stable equivalence of
Morita type between $A$ and $B$. Then, it follows from  Lemma
\ref{2.6}(2) that $(M\otimes_B-, N\otimes_A-)$ and $(N\otimes_A-,
M\otimes_B-)$ are adjoint pairs.

First, we note that $\add(Ae)=\add(Mf)$ and
$\add(N\otimes_AMf)=\add(Bf)$. In fact, this follows from the
following equalities:
$\add(Ae)=\add(M\otimes_BN\otimes_AAe)=\add(M\otimes_BBf)=\add(Mf)$,
and the fact that $\add(N\otimes_AX)=\add(N\otimes_A\add(X))$ for
any $A$-module $X$.

Thus, if a statement for the idempotent element $e$ holds true, then
it can be proved similarly for $f$, and vice versa.

Second, we shall show that the bimodules $eMf$ and $fNe$ satisfy the
conditions of a stable equivalence of Morita type between $eAe$ and
$fBf$.

$(1)$ $fNe$ is projective as an $fBf$-module and as a right
$eAe$-module, respectively.

In fact,  we have $fNe \simeq \Hom_B(Bf, {}_BNe)$ as
$fBf$-$eAe$-bimodules. Since $Ne \in \add(Bf)$ by the definition of
$f$, we see that $\Hom_{B}(Bf, Ne)$ is projective as an
$fBf$-module, that is, $fNe$ is projective as an $fBf$-module. To
see that $fNe$ is a projective right $eAe$-module, we notice that
$\add(Mf)=\add(M\otimes_BBf)=\add(M\otimes_BNe)=\add(Ae)$, here we
use the assumption $M\otimes_BNe\in \add(Ae)$.  Since $(M \otimes_B
-, N \otimes_A - )$ is an adjoint pair, it follows from
$\Hom_{B}(Bf, {}_{B}N\otimes_{A}Ae)\simeq \Hom_{A}(M\otimes_B Bf,
Ae)\simeq \Hom_A(Mf,Ae)$ that $fNe$ is projective as a right
$eAe$-module since $Mf \in \add(Ae)$. Thus (1) is proved.

$(2)$ $eMf$ is projective as an $eAe$-module and as a right
$fAf$-module, respectively. The proof of (2) is similar to that of
(1), we omit it here.

$(3)$ $eMf\otimes_{fBf}fNe\simeq eAe\oplus ePe$ as bimodules.

Indeed, by the associativity of tensor products, we have the
following isomorphisms of $eAe$-$eAe$-bimodules:
$$\begin{array}{rl} eMf\otimes_{fBf} fNe & \simeq eM\otimes_BBf\otimes_{fBf}fB\otimes_{B}Ne \\
& \simeq eM\otimes_{B}Bf\otimes_{fBf}\Hom(Bf, B)\otimes_{B}Ne\\
& \simeq eM\otimes_BBf\otimes_{fBf}\Hom(Bf, {}_BNe) \quad (\mbox{ by
Lemma \ref{2.4}}\;)\\
& \simeq eM\otimes_BNe \quad (\mbox{by Lemma \ref{2.3}}\;).
\end{array}$$
Since $M$ and $N$ define the stable equivalence of Morita type
between $A$ and $B$, we have $M\otimes_BN\simeq A\oplus P$ as
$A$-$A$-bimodules. This implies that $eMf\otimes_{fBf} fNe \simeq
eM\otimes_BNe \simeq e(A\oplus P)e\simeq eAe \oplus ePe$ as
bimodules.

$(4)$ $ePe$ is a projective $eAe$-$eAe$-bimodule.

In fact, it suffices to show that, for any indecomposable direct
summand $P'$ of the $A$-$A$-bimodule $P$, the $eAe$-$eAe$-bimodule
$eP'e$ is projective. We assume that $eP'e\neq 0$. Since
$P\in\add(A\otimes_{k}A^{op})$, there are primitive idempotent
elements  $e_1$ and $e_2$ of $A$ such that
$P'\in\add(Ae_1\otimes_ke_2A)$. Then
$_AP'e\in\add(Ae_1\otimes_ke_2Ae)\subseteq\add(Ae_1)$. This means
that $P'e$ is a direct sum of copies of $Ae_1$. Since
$P'e\in\add(Pe)\subseteq\add(Ae)$, we have $Ae_1\in\add(Ae)$.
Consequently, $eAe_1$ is a projective $eAe$-module. Now, we show
that  $e_2Ae$ is a projective right $eAe$-module. Indeed, by Lemma
\ref{2.6}(3), we have the following isomorphisms of
$A^{op}$-modules:
$eP\simeq\Hom_A(Ae,P)\simeq\Hom_A(Ae,\Hom_A(P,A))\simeq\Hom_A(P\otimes_AAe,A)\simeq\Hom_A(Pe,
A)$. This shows that $eP\in\add(eA)$ since $_APe\in \add(Ae)$. Thus
$eP'\in \add(eA)$. Since the right $A$-module $eP'$ is a direct sum
of copies of $e_2A$, it follows that $e_2A\in\add(eA)$ and $e_2Ae\in
\add(eAe)$. Consequently, $e_2Ae$ is a projective right
$eAe$-module. Hence $eAe_1\otimes_ke_2Ae$ is a projective
$eAe$-$eAe$-bimodule, and so is its direct summand $eP'e$. This
shows that $ePe$ is a projective $eAe$-$eAe$-bimodule.

$(5)$ Similarly, we can prove that $fNe\otimes eMf \simeq fBf \oplus
fQf$ as bimodules, and that the $fBf$-$fBf$-bimodule $fQf$ is
projective.

Thus, by definition, the bimodules $eMf$ and $fNe$ define a stable
equivalence of Morita type between $eAe$ and $fBf$.

Finally, the last statement of Theorem \ref{thm2} follows from
Proposition \ref{prop10} below, which emphasizes the view of
functors.

Before we give the formulation of Proposition \ref{prop10}, we
introduce here a few more notations: Set $\Lambda=\End_{eAe}(eA)$,
$R=\End_{fBf}(fN)$, $\Gamma=\End_{fBf}(fB)$,
$N'=\Hom_{fBf}((fB)_{\Gamma},fNe\otimes_{eAe}(eA)_\Lambda)$ and
$M'=\Hom_{eAe}((eA)_\Lambda,eMf\otimes_{fBf}(fB)_{\Gamma})$.

Let $\varphi:A\ra \Lambda$ be the algebra homomorphism defined by
sending $a\in A$ to $\varphi_a$, where $\varphi_a: eA\ra eA,
ex\mapsto exa$ for $x\in A$. Similarly, we define an algebra
homomorphism $\psi: B\ra \Gamma$.

Recall that, given a diagram of functors between categories:
$$\xymatrix{
{\mathcal  A}\ar^{F}[r]\ar_{H}[d] & {\mathcal B}\ar^G[d]\\
{\mathcal C}\ar^K[r] & {\mathcal D},
 }$$ we say that this diagram
is commutative if there is a natural isomorphism $\alpha: GF\ra KH.$

\begin{Prop}\label{prop10}
$(1)$ The following diagram  of functors is commutative
$$\xymatrix{
A\modcat\ar^{N\otimes_A-}[rr]\ar^{e\cdot}[d]   &&    B\modcat\ar^{M\otimes_B-}[rr]\ar^{f\cdot}[d]    &&    A\modcat\ar^{e\cdot}[d]\\
eAe\modcat\ar^{fNe\otimes_{eAe}-}[rr]         &&
fBf\modcat\ar^{eMf\otimes_{fBf}-}[rr]         && eAe\modcat.}
%\leqno(*)
$$ In particular, $_{fBf}fNe\otimes_{eAe}eA\simeq _{fBf}fN$
and $_{eAe}eMf\otimes_{fBf}fB\simeq _{eAe}eM$.

$(2)$ We have the following commutative diagram of functors
$$\xymatrix{
A\modcat\ar^{N\otimes_A-}[rr]\ar^{_\Lambda\Lambda\otimes_A-}[d] &&
B\modcat\ar^{M\otimes_B-}[rr]\ar^{_\Gamma\Gamma\otimes_B-}[d] &&
 A\modcat\ar^{_\Lambda\Lambda\otimes_A-}[d]\\
\Lambda\modcat\ar^{_\Gamma N'\otimes_\Lambda-}[rr]         &&
\Gamma\modcat\ar^{_\Lambda M'\otimes_\Gamma-}[rr] &&
\Lambda\modcat,}$$ where the right  $A$-module structure on
$\Lambda$ and the right $B$-module structure on $\Gamma$ are induced
by $\varphi$ and $\psi$, respectively. Moreover, $_\Gamma
N'_\Lambda$ and $_\Lambda M'_\Gamma$ define a stable equivalence of
Morita type between $\Lambda$ and $\Gamma$.
\end{Prop}

\textbf{Proof.} (1) To prove that the first square in (1) is
commutative, it is sufficient to show that there is a natural
transformation $\Phi: fNe\otimes_{eAe}e(-)\lra fN\otimes_A-$, which
is an isomorphism. Now we define $\Phi$ to be the composition of the
following two natural transformations: for each $X\in A$-mod,
$$\Phi_X:
fNe\otimes_{eAe}eX\lraf{\sim}fN\otimes_AAe\otimes_{eAe}eX\lraf{id_{fN}\otimes\mu}
fN\otimes_AX,$$ where $\mu: Ae\otimes_{eAe}eX\ra X$ is the
multiplication map. Clearly, we need only to show that
$id_{fN}\otimes\mu$ is a natural isomorphism, that is, for each
$_AX$, we have to show that
$$ fN\otimes_AAe\otimes_{eAe}eX\lra fN\otimes_AX$$ is an
isomorphism.

Indeed, we shall first show that if $Z\in A\modcat $ and $eZ=0$,
then $fN\otimes_A Z=0$. To prove this, we observe that $
fN\otimes_AZ\simeq\Hom_B(Bf,N\otimes_AZ)\simeq\Hom_A(_AM\otimes_B
Bf,Z), $ where the second isomorphism comes from the adjoint pair
($M\otimes_B-$,$N\otimes_A-$). Since $\add(Bf)=\add(N\otimes_AAe)$
and $Pe\in\add(Ae)$,  we have $\add(M\otimes_BBf)=\add(Ae)$. Thus
$eZ=0$ implies that $fN\otimes_A Z=0$. Next, we consider the exact
sequence $$ 0\lra \mbox{Ker}(\mu)\lra Ae\otimes_{eAe}eX\lraf{\mu}
X\lra X/AeX\lra 0.$$ Note that $e\,\mbox{Ker}(\mu)=0=e(X/AeX)$ and
that $fN_A\simeq fB\otimes_BN_A$ is a projective right $A$-module.
By applying tensor functor $fN\otimes_A-$ to the above sequence, we
deduce that
$$\begin{CD}fN\otimes_AAe\otimes_{eAe}eX@>{id_{fN}\otimes\mu}>>fN\otimes_AX\end{CD}$$
is an isomorphism. Thus we have proved the commutativity of the left
square in (1).

Similarly, we can prove that  the right square of $(\ast)$ commutes.
In particular, we see that $fNe\otimes_{eAe}eA\simeq fN_A $ as
$fBf$-$A$-bimodules, and $eMf\otimes_{fBf}fB\simeq eM_B $ as
$eAe$-$B$-bimodules.

(2) Note that  the bimodules
$\Hom_{eAe}(eA_\Lambda,eMf\otimes_{fBf}(fN)_R)$ and
$\Hom_{fBf}(fN_R,fNe\otimes_{eAe}(eA)_\Lambda)$ have been
constructed in \cite [Theorem 1.1]{LX3}, which induced a stable
equivalence of Morita type between $\Lambda$ and $R$. Since
$\add(fN)= \add(fB)$, we see that  $\Hom_{fBf}(fB,fN)$ and
$\Hom_{fBf}(fN,fB)$ induce a Morita equivalence between $R$ and
$\Gamma$. As a result,  $N'$ and $M'$ define a stable equivalence of
Morita type between $\Lambda$ and $\Gamma$. It can be checked
directly that ${_\Gamma} N'\otimes_\Lambda\Lambda_A\simeq {_\Gamma}
N'_A$ and $ {_\Lambda} M'\otimes_\Gamma\Gamma_B\simeq {_\Lambda}
M'_B$.  So, we have
$$\begin{array}{rl}
{_\Gamma} N'\otimes_\Lambda\Lambda_A &\simeq {}_\Gamma N'_A \simeq \Hom_{fBf}(fB,fN_A)\\
&\simeq  \Hom_{fBf}(fB,fB\otimes_BN_A)  \\
&\simeq{_\Gamma}\Gamma\otimes_BN_A
\end{array}$$ and
$$\begin{array}{rl}
{_\Lambda} M'\otimes_\Gamma\Gamma_B & \simeq {}_{\Lambda}M'_B = \Hom_{eAe}(eA, eMf\otimes_{fBf}fB_B)\\
                                    & \simeq\Hom_{eAe}(eA_\Lambda,eM_B) \\
                                    & \simeq\Hom_{eAe}(eA_\Lambda,eA\otimes_AM_B)\\
                                    & \simeq {_\Lambda}\Lambda\otimes_AM_B.
\end{array}$$
This implies that the diagram  in (2) is commutative. Thus, we have
proved Proposition \ref{prop10}. This also finishes the proof of
Theorem \ref{thm2}. $\square$

\medskip
{\it Remarks}. $(1)$ In Theorem \ref{thm2}, the assumption that  $M$
and $N$ do not have any projective bimodules as direct summands is
actually a very mild restriction. In fact, if $M=X'\oplus X''$ and
$N=Y'\oplus Y''$ such that $X'$ and $Y'$ have no direct summands of
projective bimodules, and that $X''$ and $Y''$ are projective
bimodules, then it follows from \cite[Lemma 4.8]{LX2} that the
bimodules $X'$ and $Y'$ also define a stable equivalence of Morita
type between $A$ and $B$. Clearly, we have $X'\otimes_BY'e\in
\add(Ae)$ and $\add(Y'e)\subseteq\add(Ne)$. Since $_AX'\otimes_BNe$
is a direct summand of $_AM\otimes_BNe$, we get
$X'\otimes_BNe\in\add(Ae)$, and
$Y'\otimes_AX'\otimes_BNe\in\add(Y'\otimes_AAe)=\add(Y'e)$. This
gives $Ne\in\add(Y'e)$. Hence $\add(Y'e)=\add(Ne)$. This means that
$M$ and $N$ in Theorem \ref{thm2} can be replaced by the bimodules
$X'$ and $Y'$.

(2) Note that $M\otimes_BNe\in \add(Ae)$ is equivalent to $Pe\in
\add(Ae)$. In Theorem \ref{thm2}, if $e$ is an idempotent element in
$A$ such that every indecomposable projective-injective $A$-module
is isomorphic to a summand of $Ae$, then $Pe\in \add(Ae)$. This
follows immediately from Lemma \ref{2.6}(3).

(3) As pointed out in \cite[Section 4]{DM}, if $e$ is an idempotent
in $A$ and if $f$ is an idempotent in $B$ such that $\add(Ae)$ and
$\add(Bf)$ are invariant under Nakayama functors, then $eAe$ and
$fBf$ are self-injective, and any stable equivalence of Morita type
between $A$ and $B$ induces a stable equivalence of Morita type
between $eAe$ and $fBf$. Note that we may recover this result from
Theorem \ref{thm2} since the idempotents $e$ and $f$ satisfy the
assumptions of Theorem \ref{thm2} by \cite[Lemma 4.1]{DM}.  In
general, however, our algebras $eAe$ and $fBf$ in Theorem \ref{thm2}
may not be self-injective.

\begin{Def}$\cite{ARSI}$
Let $A$ be an algebra. A projective $A$-module $W$ is called a
minimal Wedderburn projective module if
$\add(\nu_A(W))=\add(I_0(A)\oplus I_1(A))$, where $\nu_A$ is the
Nakayama functor of $A$ and $0\ra A\ra I_0(A)\ra I_1(A)$ is the
minimal injective copresentation of $_AA$. An idempotent element
$e\in A$ is called a minimal Wedderburn idempotent element if $Ae$
is a minimal Wedderburn projective module.
\end{Def}

Auslander  proved in \cite{ARSI} that, given $e^2 = e\in A$, the
canonical map $\rho: A\ra\End_{eAe}(eA)$ defined by right
multiplication is an isomorphism if and only if $\add(Ae)$ contains
a minimal Wedderburn projective $A$-module.

The following result shows that stable equivalences of Morita type
preserve minimal Wedderburn projective modules or minimal Wedderburn
idempotent elements.

\begin{Lem}\label{2.12}
Suppose that $A$ and $B$ are $k$-algebras such that $A$ and $B$ have
no semisimple direct summands. Assume that $_AM_B$ and $_BN_A$ do
not possess any projective bimodules as direct summands, and induce
a stable equivalence of Morita type between $A$ and $B$. Take a
minimal Wedderburn idempotent $e\in A$ and a minimal Wedderburn
idempotent $f\in B$. Then we have
$$
\add(M\otimes_BBf)=\add(Ae) \quad \text{and}\quad
\add(N\otimes_AAe)=\add(Bf).
$$
\end{Lem}
{\bf{Proof}}. We assume that $M\otimes_{B}N\simeq {A\oplus P}$ as
$A$-$A$-bimodules for some projective $A$-$A$-bimodule $P$, and
$N\otimes_{A}M \simeq {B\oplus Q}$ as $B$-$B$-bimodules for some
projective $B$-$B$-bimodule $Q$. Note that, by Lemma \ref{2.6}, the
images of the functors $_AP\otimes_A-$ and $_BQ\otimes_B-$ consist
of  projective-injective modules.

Let $0\ra A\ra I_0\ra I_1$ and $0\ra B\ra J_0\ra J_1$ be minimal
injective co-presentations of $_AA$ and $_BB$, respectively. We
claim that
$$
\add(M\otimes_B(J_0\oplus J_1))=\add(I_0\oplus I_1) \quad
\text{and}\quad \add(N\otimes_A(I_0\oplus I_1))=\add(J_0\oplus J_1).
$$
Clearly, for any $A$-module $X$, we have
$\add(N\otimes_A\add(X))=\add(N\otimes_A X)$. Since $0\ra A\ra
I_0\ra I_1$ is exact and $N_A$ is projective, it follows that
$$0\lra _BN\lra N\otimes_A I_0\lra N\otimes_A I_1$$ is exact.
Since $_BN\otimes _A DA$ is injective and $\add(_BB)=\add(_BN)$, we
see that $\add(J_0\oplus J_1)\subseteq\add(N\otimes_A(I_0\oplus
I_1))$. This implies that $\add(M\otimes_B(J_0\oplus
J_1))\subseteq\add(M\otimes_B N\otimes_A(I_0\oplus I_1))$.  Since
$P\otimes_A DA$ is projective-injective and since all indecomposable
projective-injective $A$-modules occur in $I_0$, we have
$\add(M\otimes_B N\otimes_A(I_0\oplus I_1))=\add(I_0\oplus I_1)$.
Thus, $\add(M\otimes_B(J_0\oplus J_1))\subseteq\add(I_0\oplus I_1)$.
Furthermore, it follows from the injectivity of the module
$_AM\otimes _BDB$ and $\add(_AA)=\add(_AM)$ that $\add(I_0\oplus
I_1)\subseteq\add(M\otimes_B(J_0\oplus J_1))$. Thus
$\add(M\otimes_B(J_0\oplus J_1))=\add(I_0\oplus I_1)$. Similarly, we
can prove that $\add(N\otimes_A(I_0\oplus I_1))=\add(J_0\oplus
J_1)$. Since $e\in A$ and  $f\in B$ are minimal Wedderburn
idempotents, we see that $\add(I_0\oplus I_1)=\add(\nu_A(Ae))$ and
$\add(J_0\oplus J_1)=\add(\nu_B(Bf))$. Consequently,
$\add(N\otimes_A\nu_A(Ae))= \add(\nu_B(Bf))$. It follows from
$N\otimes_A\nu_A(Ae)\simeq \nu_B(N\otimes{_A}{Ae})$ that
$\add(\nu_B(N\otimes{_A}{Ae}))= \add(\nu_B(Bf))$. Since  the
Nakayama functor $\nu_B$ is an equivalence from $B\pmodcat$ to
$B\imodcat$, we deduce that $\add(N\otimes_A Ae)=\add(Bf)$.
Similarly, we can show that $\add(M\otimes_BBf)=\add(Ae)$. $\square$

In the following we shall see that stable equivalences of Morita
type can be transfer to ``corner" algebras of Wedderburn type.

\begin{Koro}\label{th12} Suppose that $A$ and $B$ are $k$-algebras such that $A$ and $B$ have
no semisimple direct summands. Assume that $_AM_B$ and $_BN_A$ have
no projective bimodules as direct summands, and induce a stable
equivalence of Morita type between $A$ and $B$.  Let $e\in A$ and
$f\in B$ be minimal Wedderburn idempotents. Then $eMf$ and $fNe$
define a stable equivalence of Morita type between $eAe$ and $fBf$
such that $fNe\otimes_{eAe}eA\simeq fN$ and
$eMf\otimes_{fBf}fB\simeq eM$ as bimodules.
\end{Koro}

{\bf Proof}. By Lemma \ref{2.12}, we see that the idempotents $e$
and $f$ satisfy the assumptions in Theorem \ref{thm2}. Then
Corollary \ref{th12} follows from  the first part of Theorem
\ref{thm2} together with Proposition \ref{prop10}. $\square$

As a corollary of Corollary \ref{th12}, we get the following result.

\begin{Koro}\label{coro12} Assume that $A$ and $B$ are $k$-algebras
without semisimple direct summands. Let $_AX$ be a
generator-cogenerator for $A\emph{-mod}$, and let $_BY$ be a
generator-cogenerator for $B\emph{-mod}$.  If $\End_A(X)$ and
$\End_B(Y)$ are stably equivalent of Morita type, then there exist
bimodules $_AM_B$ and $_BN_A$ which define a stable equivalence of
Morita type between $A$ and $B$ such that $\add(_AM\otimes_B
Y)=\add(_AX)$ and $\add(_BN\otimes_A X)=\add(_BY)$.
\end{Koro}

{\bf Proof}. Set $R= \End_A(X)$ and $S= \End_B(Y)$. First, we show
that if $A$ does not have any semisimple direct summands, then nor
does $R$.

Suppose contrarily that $R$ has a semisimple direct summand. Then
$R$ must have a simple projective-injective module $W$. Since each
indecomposable projective-injective $R$-module is isomorphic to a
direct summand of $\Hom_A(X,DA)$, there exists an indecomposable
injective $A$-module $I$ such that $W\simeq\Hom_A(X,I)$. Let $_AS$
be the socle of $_AI$. Then $\Hom_A(X,S)$ can be embedded into the
simple $R$-module $\Hom_A(X,I)$, and therefore
$\Hom_A(X,S)\simeq\Hom_A(X,I)\simeq W$ as $R$-modules. Since
$A\in\add(X)$, we infer that $S\simeq I$. Let $_AP$ be the
projective cover of $_AS$. Then it follows from $\Hom_R(\Hom_A(X,P),
\Hom_A(X,S)) \simeq \Hom_A(P,S)\neq 0$ that there is a non-zero
homomorphism from $\Hom_A(X,P)$ to the simple projective $R$-module
$\Hom_A(X,S)$, which means that $\Hom_A(X,P)\simeq\Hom_A(X,S)$.
Consequently, we get $P\simeq S\simeq I$. Thus $A$ has a simple
projective-injective module, and therefore it has a semisimple
direct summand, which is a contradiction. This shows that $R$ has no
semisimple direct summands. Similarly, we can prove that $S$ has no
semisimple direct summands.

Note that, if $X$ is a generator-cogenerator for $A$-mod, then
$\Hom_A(X,A)$ is a minimal Wedderburn projective $R$-module.
Similarly, $\Hom_B(Y,B)$ is a minimal Wedderburn projective
$S$-module. Clearly, $\End_{R}(\Hom_A(X,A))\simeq A$ and
$\End_S(\Hom_B(Y,B))\simeq B$. Note that neither $R$ nor $S$ has
semisimple direct summands. Then, by Corollary \ref{th12}, there
exist bimodules $_AM_B$ and $_BN_A$ which define a stable
equivalence of Morita type between $A$ and $B$. Note that
$\Hom{_R}(\Hom_A(X,A), R)\simeq {}_AX$ and $\Hom{_S}(\Hom_B(Y,B),
S)\simeq {}_BY$.  It follows  from Corollary \ref{th12} that
$\add(_AM\otimes_B Y)=\add(_AX)$ and $\add(_BN\otimes_A
X)=\add(_BY)$. $\square$

\smallskip
Combining Corollary \ref{coro12} with \cite[Theorem 1.1]{LX3}, we
have the following result on Auslander algebras.

\begin {Koro}  Let $A$ and $B$ be representation-finite $k$-algebras.
Suppose that $A$ and $B$ have no semisimple direct summands. Let
$\Lambda$ and $\Gamma$ be the corresponding Auslander algebras of
$A$ and $B$, respectively. Then $\Lambda$ and $\Gamma$ are stably
equivalent of Morita type if and only if so are $A$ and $B$.
\label{cor3}
\end{Koro}

\medskip
For an algebra $A$, we denote by $[A]$ the class of all algebras $B$
such that there is a stable equivalence of Morita type between $B$
and $A$. It is known that $[A]=[A\times S]$ for any separable
algebra $S$. Note that, if $k$ is a perfect field, then the class of
all semisimple $k$-algebras is the same as that of all separable
$k$-algebras.

The following result establishes a one-to-one correspondence, up to
stable equivalence of Morita type, between representation-finite
algebras  and Auslander algebras. This is an immediate consequence
of Corollary \ref{cor3}.

\begin {Koro}\label{coro13}
Suppose that $k$ is a perfect field. Let $\cal F$ be the set of
equivalence classes $[A]$ of representation-finite $k$-algebras $A$
under stable equivalences of Morita type, and let $\cal A$ be the
set of equivalence classes $[\Lambda]$ of Auslander $k$-algebras
$\Lambda$ under stable equivalences of Morita type. Then there is a
one-to-one correspondence between $\cal F$ and $\cal A$.
\end{Koro}

Finally, we remark that Corollary \ref{cor3} is not true for derived
equivalences. Nevertheless, it was shown in \cite{HAU} that if two
representation-finite, self-injective algebras $A$ and $B$ are
derived-equivalent then so are their Auslander algebras. The
converse of this statement is open. For further information on
constructing derived equivalences, we refer the reader to the
current papers \cite{hx2, HAU}.

\section{Stable equivalences of Morita type based on self-injective algebras\label{selfinjective}}

Of particular interest are stable equivalences of Morita type
between self-injective algebras or between those related to
self-injective algebras. Since derived equivalences between
self-injective algebras imply stable equivalences of Morita type by
a result of Rickard \cite{R3}, this makes stable equivalences of
Morita type closely  related to the Brou\'e abelian defect group
conjecture which essentially predicates a derived equivalence
between two block algebras \cite{B}, and thus also a stable
equivalence of Morita type between them.

In this section, we will apply Theorem \ref{thm1} and Theorem
\ref{thm2} to self-injective algebras. It turns out that the
existence of a stable equivalence of Morita type between
$\Phi$-Auslander-Yoneda algebras of generators for one finite
admissible set $\Phi$ implies the one for all finite admissible
sets.

Throughout this section, we fix a finite admissible subset $\Phi$ of
$\mathbb{N}$, and assume that  $A$ and $B$ are indecomposable,
non-simple, self-injective algebras. Let $X$ be a generator for
$A$-mod with a decomposition $X:= A\oplus \displaystyle\bz_{1\leq
i\leq n} X_i$ , where $X_i$ is indecomposable and non-projective
such that $X_i\ncong X_t$ for $1\leq i\neq t\leq n$, and let $Y$ be
a generator for $B$-mod with a decomposition $Y:= B\sz
\displaystyle\bz_{1\leq j\leq m} Y_i$, where
 $Y_j$ is indecomposable and non-projective such that $Y_j\ncong
Y_s$ for $ 1\leq j\neq s\leq m$.

\begin{Lem} \label{lem 2.13}
$(1)$ The full subcategory of\, $\E_A^{\Phi}(X)\modcat$ consisting
of projective-injective  $\E_A^{\Phi}(X)$-modules is equal to
$\add(\E_A^{\Phi}(X, A))$. Particularly, if
$\E_A^{\Phi}(X)\neq\End_A(X)$, then $\dm(\E_A^{\Phi}(X))= 0$.

$(2)$ $\E_A^{\Phi}(X)$ has no semisimple direct summands.
\end{Lem}

{\bf Proof}. (1) For convenience, we set $\Lambda_0= \End_A(X)$ and
$\Lambda= \E_A^{\Phi}(X)$.  Since $A$ is self-injective, it follows
from \cite[\,Lemma\,3.5]{HAU} that $\nu{_\Lambda}(\E_A^{\Phi}(X,
A))\simeq \E_A(X,\nu_AA)\simeq \E_A^{\Phi}(X, DA)\in \add\big(
\E_A^{\Phi}(X, A)\big)$. Consequently, $\E_A^{\Phi}(X, A)$ is a
projective-injective $\Lambda$-module. We claim that, up to
isomorphism, each indecomposable projective-injective
$\Lambda$-module is a direct summand of $\E_A^{\Phi}(X, A)$. To
prove this claim, it suffices to show that $\E_A^{\Phi}(X, X_i)$ is
not injective for all $1\leq i\leq n$. We denote $\E_A^{\Phi}(X,
X_i)$ by $\widetilde{X_i}$ for abbreviation.

First, we observe that $\rad(\Lambda)=
\rad(\Lambda_{0})\oplus\Lambda_{+}$, where $\Lambda_{+} =
\displaystyle\bz_{0\neq i\in \Phi}\Lambda_i$ with
$\Lambda_i=\Ext^i_A(X,X)=\Hom_{\Db{A}}(X,X[i])$. Since each summand
$\Hom_{\Db{A}}(X,X[j])$ of $\widetilde{X_i}$ is a $\Lambda_0$-module
and since the socle of $\widetilde{X_i}$ is the set of all elements
$x$ in $\widetilde{X_i}$ such that $\rad(\Lambda)x=0$, we see that
the socle of $\widetilde{X_i}$ contains $\displaystyle\bz_{j\in
\Phi}\{x\in \soc{_{\Lambda_0}}(\Ext_A^j(X, X_i))\mid \Lambda_{+}x=
0\}$. By an argument of graded modules, we can even see that
$\soc_{\Lambda}(\widetilde{X_i}) = \displaystyle\bz_{j\in
\Phi}\{x\in \soc{_{\Lambda_0}}(\Ext_A^j(X, X_i))\mid \Lambda_{+}x=
0\}$.

Next, we shall show that $\widetilde{X_m}$ is not injective for
$1\leq m\leq n$. Indeed, let $f: X_m\to I$ be an injective envelope
of $X_m$ with $I$ an injective $A$-module. Then $f{_\ast}: \Hom_A(X,
X_m)\to\Hom_A(X, I)$ is an injective envelop of the
$\Lambda_0$-module $\Hom_A(X,X_m)$ in $\Lambda_{0}\modcat$. Now, we
consider the following two cases:

(a) If $\widetilde{X_m} = \Hom_A(X, X_m)$, then $\widetilde{X_m}$ is
annihilated by $\Lambda_{+}$. Since $X_m$ is not injective in
$A\modcat$, we conclude that $\Hom_A(X, X_m)$ is not an injective
$\Lambda_0$-module, which implies that $\widetilde{X_m}$ is not
injective as a $\Lambda$-module.

(b) If $\widetilde{X_m}\neq \Hom_A(X, X_m)$, then there is a
positive integer $ t\in\Phi$ such that $\Ext_A^t(X, X_m)\neq 0$. We
may assume that $t$ is the maximal number in $\Phi$ with this
property, that is, $\Ext_A^s(X, X_m)= 0$ for any $s\in\Phi$ with
$t<s$. It follows that $\Lambda_{+}\Ext_A^t(X, X_m)= 0$, which
implies that $0\neq \,\soc{_{\Lambda_0}}(\Ext_A^t(X, X_m))\subseteq
\soc{_\Lambda}(\widetilde{X_m})$.

Now we consider $\soc{_{\Lambda_0}}(\Hom_A(X, X_m))$. Since
$f{_\ast}$ is an injective envelop in $\Lambda_{0}\modcat$, we know
that $ \soc{_{\Lambda_0}}(\Hom_A(X, X_m))\simeq
\soc{_{\Lambda_0}}(\Hom_A(X, I))$. Since
$\nu{_{\Lambda_0}}(\Hom_A(X, A))\in \add\big(\Hom_A(X, A)\big)$ and
$I\in\add(_AA)$, we see that $\Hom{_{\Lambda_0}}(\Hom_A(X, X_i),
\soc{_{\Lambda_0}}(\Hom_A(X, I))) = 0$ for $1\leq i\leq n$. If $e$
is the idempotent in $\Lambda_{0}$ corresponding to the direct
summand $A$ of $X$, then $e\,\soc{_{\Lambda_0}}(\Hom_A(X,
I))=\soc{_{\Lambda_0}}(\Hom_A(X, I))$. Consequently,
$e\,\soc{_{\Lambda_0}}(\Hom_A(X, X_m))=\soc{_{\Lambda_0}}(\Hom_A(X,
X_m))$, that is, $eg=g$  whenever $g\in\soc{_{\Lambda_0}}(\Hom_A(X,
X_m))$, that is, $g$ factorizes through the regular module $_AA$,
say $g=g_1g_2$ with $g_1: X\ra {}_AA$ and $g_2:{}_AA\ra X_m$. Thus,
for any element $x\in \Hom_{\Db{A}}(X,X[i])$ with $0\ne i\in \Phi$,
we have $x\cdot g=x(g_1[i]\;g_2[i])=(x\;g_1[i])g_2[i]=0\,(g_2[i])=0$
since $A$ is self-injective. Thus
$\Lambda_+\soc_{\Lambda_0}(\Hom_A(X,X_m))=0$. This implies that
$\soc_{\Lambda_0}(\Hom_A(X,X_m))\subseteq
\soc_{\Lambda}(\widetilde{X_m})$. Thus we have shown that the
$\Lambda$-submodule $\soc_{\Lambda_0}(\Hom_A(X,X_m))\oplus
\soc_{\Lambda_0}(\Ext^t_A(X,X_m))$ of $\widetilde{X_m}$ is contained
in the socle of $\widetilde{X_m}$. This implies that
$\widetilde{X_m}$ is not injective since its socle is not simple.

Thus $\add(\E_A^{\Phi}(X, A))$ is just the full subcategory of\,
$\E_A^{\Phi}(X)\modcat$ consisting of projective-injective modules.

Finally, we consider the dominant dimension of
$\dm(\E_A^{\Phi}(X))$. Suppose $\E_A^{\Phi}(X)\neq\End_A(X)$. Since
$A$ is self-injective, we have $\E_A^{\Phi}(X, A)=\Hom_A(X, A)$.  It
follows that $\E_A^{\Phi}(X, A)$ is annihilated by $\Lambda_{+}$,
but not by $\Lambda$. Hence $\Lambda$ cannot be cogenerated by
$\E_A^{\Phi}(X, A)$. This implies that $\dm(\E_A^{\Phi}(X))= 0$. We
finish the proof.

(2) Contrarily, we suppose that the algebra $\E_A^{\Phi}(X)$ has a
semisimple direct summand. Then $\E_A^{\Phi}(X)$ has a simple
projective-injective module $S$. According to $(1)$, we know that
$S$ must be a simple projective-injective $\End_A(X)$-module. Then
it follows from the first part of the proof of Corollary
\ref{coro12} that $A$ has a semisimple direct summand. Clearly, this
is contrary to our initial assumption that $A$ is indecomposable and
non-simple. Thus $\E_A^{\Phi}(X)$ has no  semisimple direct
summands. $\square$

\begin{Theo}\label{coro14}
If the algebras $\E_A^{\Phi}(X)$ and $\E_B^{\Phi}(Y)$ are stably
equivalent of Morita type, then $n=m$ and  there are bimodules
$_AM_B$ and $_BN_A$ which define a stable equivalence of Morita type
between $A$ and $B$ such that, up to the ordering of indices,
\,$_AM\otimes_B Y_i\simeq {X_i\oplus P_i}$ as $A$-modules, where
$_AP_i$ is projective for all $i$ with $ 1\leq i\leq n$. Moreover,
for any finite admissible subset $\Psi$ of $\mathbb{N}$, there is a
stable equivalence of Morita type between $\E_A^{\Psi}(X)$ and
$\E_B^{\Psi}(Y)$.
\end{Theo}

{\bf Proof}. For convenience, we set $\Lambda_0= \End_A(X),
\,\Lambda= \E_A^{\Phi}(X),\, \Gamma_0= \End_B(Y)$ and $\Gamma=
\E_B^{\Phi}(Y)$. By Lemma \ref{lem 2.13}, the algebras $\Lambda$ and
$\Gamma$ have no semisimple direct summands.  Let $e$ be the
idempotent in $\Lambda_{0}$ corresponding to the direct summand $A$
of $X$, and let $f$ be the idempotent in $\Gamma_{0}$ corresponding
to the direct summand $B$ of $Y$. Note that $\Lambda e\simeq
\E_A^{\Phi}(X, A)$ as $\Lambda$-modules  and $\Gamma f \simeq
\E_B^{\Phi}(Y, B)$ as $\Gamma$-modules. Clearly, $e\Lambda e\simeq
A$ and $f\, \Gamma f \simeq B$ as algebras. Moreover, we see that
$e\Lambda \simeq X$ as $A$-modules, and $f\, \Gamma \simeq Y$ as
$B$-modules. Suppose that a stable equivalences of Morita type
between $\Lambda$ and $\Gamma$ is given. By Corollary \ref{2.8} and
Lemma \ref{lem 2.13}, we know that the idempotent $e$ in $\Lambda$
and the idempotent $f$ in $\Gamma$ satisfy the conditions in Theorem
\ref{thm2}. It follows from Theorem \ref{thm2} and Proposition
\ref{prop10}(1) that there are bimodules $_AM_B$ and $_BN_A$ which
define a stable equivalence of Morita type between $A$ and $B$ such
that $\add(M\otimes_B Y)= \add(X)$. By the given decompositions of
$X$ and $Y$, we conclude that $n=m$ and, up to the ordering of
direct summands, we may assume that \,$_AM\otimes_B Y_i\simeq
{X_i\oplus P_i}$ as $A$-modules, where $_AP_i$ is projective for all
$i$ with $ 1\leq i\leq n$. Now, the last statement in this corollary
follows immediately from Theorem \ref{thm1}. Thus the proof is
completed. $\square$

\medskip
Usually, it is difficult to decide whether an algebra is not stably
equivalent of Morita type to another algebra. The next corollary,
however, gives a sufficient condition to assert when two algebras
are not stably equivalent of Morita type.

\begin{Koro}\label{coro15}
Let $n$ be a non-negative integer. Let $W$ be an indecomposable
non-projective $A$-module. Suppose that $\Omega_A^{s}(W)\not\simeq
W$ for any non-zero integer $s$. Set $W_n= \bz _{0\leq i\leq n}
\Omega_A^{i}(W)$. Then, for any finite admissible subset $\Psi$ of
$\mathbb{N}$,  the algebras $\E_A^{\Psi}(A \oplus W_n \oplus
\Omega_A^{l}(W))$ and $\E_A^{\Psi}(A \oplus W_n\oplus
\Omega_A^{m}(W))$ are not stably equivalent of Morita type whenever
$m $ and $l$ belong to $\mathbb{N}$ with $n< m< l$.
\end{Koro}
{\bf{Proof}}. Suppose that there is  a finite admissible subset
$\Psi$ of $\mathbb{N}$ such that $\E_A^{\Psi}(A \oplus W_n \oplus
\Omega_A^{m}(W))$ and $\E_A^{\Psi}(A \oplus W_n\oplus
\Omega_A^{l}(W))$ are  stably equivalent of Morita type for some
fixed $l, m\in\mathbb{N}$ with $n< m< l$. Set $\Phi_1=\{0,1,\cdots,
n\}\cup\{l\}$ and $\Phi_2=\{0,1,\cdots,n\}\cup\{ m\}$. Then, by
Theorem \ref{coro14}, we know that there exist bimodules $_AM_A$ and
$_AN_A$ which define a stable equivalence of Morita type between $A$
and itself, and  that there is a bijection $\sigma: \Phi_1\to
\Phi_2$ such that $M \otimes_A \Omega_{A}^{j}(W) \simeq
\Omega{_{A}^{\sigma(j)}}(W) \oplus P_j$ as $A$-modules, where $P_j$
is  projective  for each  $j\in \Phi_1$. In particular, we have $M
\otimes_A W \simeq \Omega{_{A}^{\sigma(0)}}(W) \oplus P_0$. Since
$M$ is  projective as a one-sided module, we know that
$M\otimes_A\Omega_{A}^{l}(W) \simeq \Omega{_{A}^{\sigma(0)+l}}(W)
\oplus P_l'$ with $P_l' \in \add(_AA)$. Note that  $M \otimes_A
\Omega_{A}^{l}(W) \simeq \Omega{_{A}^{\sigma(l)}}(W) \oplus P_l$. It
follows that $\Omega{_{A}^{\sigma(0)+l}}(W) \simeq
\Omega{_{A}^{\sigma(l)}}(W) $. Consequently, we have $\sigma(l)=
\sigma(0)+l \geq l$ since $W$ is not $\Omega$-periodic. Hence $l\le
\sigma(l)\leq m < l$, a contradiction. This shows that
$\E_A^{\Psi}(A \oplus W_n \oplus \Omega_A^{m}(W))$ and
$\E_A^{\Psi}(A \oplus W_n\oplus \Omega_A^{l}(W))$ cannot be stably
equivalent of Morita type whenever $l$ and $m$ $\in\mathbb{N}$ with
$n< m < l$. $\square$

\medskip
This corollary will be used in the next section.

\section{A family of derived-equivalent algebras: application
to Liu-Schulz algebras\label{family}}

In this section, we shall apply our results in the previous sections
to solve the following problem on derived equivalences and stable
equivalences of Morita type:

\medskip
{\bf Problem.} Is there any infinite series of finite-dimensional
$k$-algebras  such that they have the same dimension and are all
derived-equivalent, but not stably equivalent of Morita type ?

\medskip
This problem was originally asked by Thorsten Holm at a workshop in
Goslar, Germany.

Recall that Liu and Schulz in \cite{Ls} constructed a local
symmetric $k$-algebra $A$ of dimension 8 and an indecomposable
$A$-module $M$ such that all the syzygy modules $\Omega_A^n (M)$
with $n\in \mathbb{Z}$ are 4-dimensional and pairwise
non-isomorphic.  This algebra $A$ depends on a non-zero parameter
$q\in k$, which is not a root of unity, and has an infinite
DTr-orbit in which each module has the same dimension. A thorough
investigation of Auslander-Reiten components of this algebra was
carried out by Ringel in \cite{R}. Based on this symmetric algebra
and a recent result in \cite{hx2} together with the results in the
previous sections, we shall construct an infinite family of
algebras, which provides a positive solution to the above problem.

From now on, we fix a non-zero element $q$ in the field $k$, and
assume that $q$ is not a root of unity. The $8$-dimensional
$k$-algebra $A$ defined by Liu-Schulz is an associative algebra
(with identity) over $k$ with

the generators: $x_0, x_1, x_2$, and

the relations: $ x_i^2= 0, \quad \mbox{and}\quad x_{i+1}x_i+q
x_ix_{i+1}= 0 \quad \mbox{for}\quad  i= 0, 1, 2. $

\noindent Here, and in what follows, the subscript is modulo $3$.

Let $n$ be a fixed natural number, and let $\Phi=\{0\}$ or $\{0,
1\}$. For $j\in \mathbb{Z}\,$, set $u_j: =x_2+q^jx_1 $, $I_j :=
Au_j$, $J_j:=u_jA$, $I:= \bz\limits_{i=0}^{n}I_i$ and
$\Lambda_j^{\Phi}:= \E_A^{\Phi}(A \sz I \sz I_j )$.

With these notations in mind, the main result in this section can be
stated as follows:

\begin{Theo} For any  $m \geq n+4$,  we have

\medskip
$(1)$ $\dk{(\Lambda_m^{\Phi})}= \dk{(\Lambda_{m+1}^{\Phi})}$.

\medskip
$(2)$ $\gd(\Lambda_m^{\Phi})= \infty$.

\medskip
$(3)$ $\dm(\Lambda_m^{\Phi})=\left\{\begin{array}{ll} 2 & \mbox{if } \Phi=\{0\},\\
0 & \mbox{if } \Phi=\{0,1 \}.\end{array} \right.$

\medskip
$(4)$  $\Lambda_m^{\Phi}$ and $\Lambda_{m+1}^{\Phi}$ are
derived-equivalent.

\medskip
$(5)$ If $l>m$, then $\Lambda_l^{\Phi}$ and $\Lambda_m^{\Phi}$ are
not stably equivalent of Morita type. \label{th1}
\end{Theo}

An immediate consequence of Theorem \ref{th1} is the following
corollary, which solves the above mentioned problem positively.

\begin{Koro}
There exists an infinite series of finite-dimensional $k$-algebras
$A_i,\,i\in\mathbb{N}$, such that

$(1)$ $\dk(A_i)= \dk(A_{i+1})$ for all $i\in {\mathbb N}$,

$(2)$ all $A_i$ have the same global and dominant dimensions,

$(3)$ all $A_i$ are derived-equivalent, and

$(4)$ $A_i$ and $A_j$ are not stably equivalent of Morita type for
$i\neq j$. \label{coro16}
\end{Koro}

The proof of Theorem \ref{th1} will cover the rest of this section.
Let us first introduce a few more notations and conventions.

Let $B$ be an algebra and $S$ a subset of $B$. Set $R(S):=\{b\in
B\mid sb=0 \;\mbox{for all}\; s\in S\}$ for the right annihilator of
$S$ in $B$, and $L(S):=\{b\in B\mid bs=0 \;\mbox{for all}\; s\in
S\}$ for the left annihilator of $S$ in $B$. In case $x\in B$, we
write $R(x)$ and $L(x)$ for $R(\{x\})$ and $L(\{x\})$, respectively.
For $y,z\in B$, we set $B(y, z):=\{b\in B~|~L(y)bz=0\}.$ Note that
$L(S)$ and $R(S)$ are left and right ideals in $B$, respectively.

Let $V$ be a $k$-vector space with $y_i\in V$ for $1\leq i\leq n\in
\mathbb{N}$. We denote by $<y_1,\ldots,y_n>$ the $k$-subspace of $V$
generated by all $y_i$.

The following result is useful for our calculations, it may be of
its own interest in describing the endomorphism rings of direct sums
of cyclic left ideals.

\begin{Lem}
Let $B$ be a $k$-algebra, and let $x$ , $y$  and $z$ be  elements in
$B$. Then the following statements hold:

$(1)$ There is an isomorphism of $k$-vector spaces:
$$
\varphi_{x, y}: \Hom_B(Bx, By) \lraf{\sim} R(L(x))\cap By,
$$
which sends $f$ to $f(x)$ for $f\in\Hom_B(Bx, By)$.

$(2)$ There is an isomorphism of $k$-vector spaces:
$$
\theta_{x, y}: \Hom_B(Bx, By) \lraf{\sim}B(x,y)/L(y),
$$
which sends $h$  to $\overline{b}$ for $h\in\Hom_B(Bx, By)$, where
$b\in B$ such that $h(x)=by$ and $\ol{b}$ stands for the coset
$b+L(y)$.

$(3)$ The maps $\theta_{x, x}$ and $\theta_{y, y}$ are isomorphisms
of algebras.

$(4)$ The map $\theta_{x, y}$ satisfies the following identity:
$$\theta_{x,y}(agc)=\theta_{x,x}(a)\theta_{x,y}(g)\theta_{y,y}(c)$$
for $a\in\End_A(Bx), g\in\Hom_B(Bx,By)$ and $c\in\End_B(By)$.

$(5)$ The following diagram is commutative:
$$
\xymatrix{ \Hom_B(Bx,By)\otimes_{\End_B(By)}\Hom_B(By,Bz)
\ar[d]^-{\wr}_-{\theta_{x,y}\otimes\theta_{y,z}}\ar[rr]^-{\triangle_{x,y,z}}
&&\Hom_B(Bx,Bz)\ar[d]^-{\wr}_-{\theta_{x,z}}\\
(B(x,y)/L(y))\otimes_{B(y,y)/L(y)}(B(y,z)/L(z))\ar[rr]^-{\nabla_{x,y,z}}
&&B(x,z)/L(z),}
$$
where $\triangle_{x,y,z}$ is the composition map, and
$\nabla_{x,y,z}$ is the multiplication map.

$(6)$ Let $n$ be a positive integer, and let $x_i$ be elements in
$B$ for $1\leq i\leq n$. We define
$${\rm{M}}_{B}(x_1,x_2,\cdots,x_n):=\{(\;\ol{b_{i,j}}\;){_{1\leq i,j\leq n}}~|~\ol{b_{i,j}}\in
B(x_i,x_j)/L(x_j)~\mbox{\,for\, all \,} 1\leq i,j\leq n~\}.$$ Then
${\rm{M}}_{B}(x_1,x_2,\cdots,x_n)$ becomes an associative
$k$-algebra with the usual matrix addition and  multiplication which
is given by $\nabla_{x_r,x_s, x_t}$ for $1\leq r,s,t\leq n$.
Moreover, there is a natural algebra isomorphism $\theta:
\End_B(\bz_{{0\leq i\leq n}}B{x_i})\lra
{\rm{M}}_{B}(x_1,x_2,\cdots,x_n)$, which is induced by
$\theta_{x_r,x_s}$ for $1\leq r,s\leq n$.

\label{2}
\end{Lem}

{\bf Proof}.$(1)$ Let $f\in\Hom_B(Bx, By)$. Since $f$ is a
homomorphism of $B$-modules, we know $b(xf)=0$ whenever $b\in B $
and $bx=0$. This implies that $xf\in R(L(x))\cap By$. Thus the map
$\varphi_{x,y}$ is well-defined. It is not hard to check that
$\varphi_{x,y}$ is an isomorphism of $k$-vector spaces.

$(2)$ For $x\in B$, we denote by $\rho_x$ the right multiplication
map from $B$ to itself, defined by $ b\mapsto bx$ for $b\in B$. Then
there is a canonical exact sequence of $B$-modules: $\delta_{x}:
0\ra L(x)\lraf{\lambda_{x}} B \lraf{\pi_{x}} Bx \ra 0$, where
$\lambda_x$ is the inclusion, and $\pi_x$ is the canonical
multiplication of $x$. By the definition of $B(x,y)$, an element $w$
belongs to $B(x,y)$ if and only if $\lambda_x\rho_w\pi_y=0$, or
equivalently, if and only if there is a unique
$\alpha\in\Hom_B(Bx,By)$ such that $\rho_w\pi_y=\pi_x\alpha$.
Clearly, $w\in L(y)$ if and only if $\rho_w\pi_y=0$. So, we have
$L(y)\subseteq B(x,y)$.

First, we show that $\theta_{x,y}$ is well-defined. In fact, if
$f\in\Hom_B(Bx,By)$, then there is an element $b\in B$, which may
not be unique, such that the following diagram of left $B$-modules
commutes:
$$
\xymatrix{
0\ar[r]&L(x)\ar@{-->}_{\rho_b'}[d]\ar[r]^-{\lambda_x}&B\ar[r]^-{\pi_x}\ar@{-->}_{\rho_b}[d]
&Bx\ar[d]^-{f}\ar[r]&0\\
0\ar[r]&L(y)\ar[r]^-{\lambda_y}&B\ar[r]^-{\pi_y}&By\ar[r]& 0,}
$$
where $\rho_b'$ is the restriction of $\rho_b$ to $L(x)$. Hence
$b\in B(x,y)$. If there is another $d$ in $B$ also making the above
diagram commutative, then $(\rho_b-\rho_d)\pi_y=0$, and therefore
$\rho_b-\rho_d$ factorizes through $L(y)$. This implies that $b-d\in
L(y)$ and $\ol{b}=\ol{d}$ in $B(x,y)/L(y)$. Thus $\theta_{x,y}$ is
well-defined.

Next, we shall prove that $\theta_{x,y}$ is an isomorphism of
$k$-vector spaces. Indeed, if $\theta_{x,y}(f)=\ol{b}=0$ for some
map $f\in \Hom_B(Bx,By)$, then $b\in L(y)$ and $\pi_xf = \rho_b\pi_y
= 0$. Since $\pi_x$ is surjective, we get $f=0$. Thus $\theta_{x,y}$
is injective. That $\theta_{x.y}$ is surjective follows from the
equivalent definitions of $B(x,y)$ discussed above.

(3) to (5) follow from the above proof of (2).

(6) is a consequence of (2) to (5). $\square$

\medskip
Recall that, for $i\in \mathbb{Z}$, we have defined
$u_i:=x_2+q^ix_1, I_i:=A u_i $ and $J_i:=u_i A$. In the following
lemma, we display a few properties about the Liu-Schulz algebra $A$.

\begin{Lem}\emph{\cite{Ls, R}}
$(1)$ The Liu-Schulz algebra $A$ is an $\mathbb N$-graded algebra,
namely, $A=\bz_{i\geq 0}A_i$ with
$$
A_0=k,\quad A_1=<x_0, x_1, x_2>,\quad A_2=<x_0x_1, x_1x_2,
x_2x_0>,\quad A_3=<x_0x_1x_2>, \quad \mbox{and}\quad
A_i=0\quad\text{for all}\quad  i\geq 4.$$ Moreover, $A_2$ is
contained in the center of $A$. In particular,
$x_0x_1x_2=x_1x_2x_0=x_2x_0x_1$ in $A$.

$(2)$ $A$ is an $8$-dimensional symmetric $k$-algebra.

$(3)$ $\dk(I_j)=\dk(J_j)=4$ for all $j\in \mathbb{Z}$.

$(4)$ $\Omega_A(I_j)=I_{j+1}$ and $\Omega_A(J_{j+1})=J_j$ for all
$j\in \mathbb{Z}$.

$(5)$ The $A$-modules $I_j$ (respectively, $A^{op}$-modules $J_j$)
are pairwise non-isomorphic  for all $j\in \mathbb{Z}$.
\label{lem2'}
\end{Lem}

In the next lemma, we calculate dimensions of homomorphism groups
related to the modules $I_i$ and $J_i$.

\begin{Lem} Let $i$ and $j$ be integers. Then

$(1)$ $I_j$ has a basis $\{x_2+q^jx_1, x_2x_0-q^{j-1}x_0x_1, x_1x_2,
x_0x_1x_2\}$, and
    $J_j$ has a basis $\{x_2+q^jx_1, x_2x_0-q^{j+1}x_0x_1, x_1x_2, x_0x_1x_2\}$.

$(2)$ $L(u_j)=I_{j+1}$, $R(u_{j+1})=J_j$.

$(3)$ $J_j\simeq \Hom_A(I_j, A)$.

$(4)$ As $k$-vector spaces, $\Hom_A(I_j, I_i)\simeq J_j\cap I_i=
\left\{\begin{array}{ll} <x_2+q^jx_1, x_1x_2, x_0x_1x_2> & \mbox{if } j=i,\\
<x_2x_0-q^{j+1}x_0x_1, x_1x_2, x_0x_1x_2> & \mbox{if } j=i-2,\\
<x_1x_2, x_0x_1x_2> & \mbox{otherwise}.\end{array} \right.$

In particular, $\dk\Hom_A(I_j, I_i)=
\left\{\begin{array}{ll} 3 & \mbox{if } j=i \mbox{ or } i-2,\\
2 & \mbox{otherwise}.\end{array} \right.$

$(5)$ $\dk\Ext_A^1(I_j,I_i)=
\left\{\begin{array}{ll} 1 & \mbox{if } j\leq i\leq j+3,\\
0 & \mbox{otherwise}.\end{array} \right.$

$(6)$ $A(1,u_i)=A$ and $A(u_i,1)=J_i$.

$(7)$ $A({u_j,u_i})=
\left\{\begin{array}{ll} <1,x_1,x_2,x_0x_1,x_1x_2,x_2x_0,x_0x_1x_2> & \mbox{if } j=i,\\
<x_1,x_2,x_0x_1,x_1x_2,x_2x_0,x_0x_1x_2> & \mbox{if } j=i-2,\\
<x_0,x_1,x_2,x_0x_1,x_1x_2,x_2x_0,x_0x_1x_2> &
\mbox{otherwise}.\end{array} \right.$
\label{3}
\end{Lem}

{\bf Proof}.  (1) and (2). By definition, $I_j=Au_j$. One can check
directly that
\begin{eqnarray*}
x_0u_j=(-q)(x_2x_0-q^{j-1}x_0x_1), & x_2u_j=-q^{j+1}x_1x_2, &
x_1u_j=x_1x_2,
\\
x_1x_2u_j= x_0x_1x_2u_j=0,& x_0x_1u_j=x_0x_1x_2, &
x_2x_0u_k=q^jx_0x_1x_2.
\end{eqnarray*}
This implies that  $ I_j=<x_2+q^jx_1, x_2x_0-q^{j-1}x_0x_1, x_1x_2,
x_0x_1x_2>$. Note that $0\ra L(u_j)\ra A\ra Au_j\ra 0$ is an exact
sequence of $A$-modules.  Since
$u_{j+1}u_j=(x_2+q^{j+1}x_1)(x_2+q^jx_1)=0$, we have $
I_{j+1}\subseteq L(u_j)$. In addition,  $\dk I_{j+1}=\dk L(u_j)=4$.
It follows that $L(u_j)=I_{j+1}$. Similarly, we can prove the
corresponding statements in (1) and (2) for $J_j$.

(3) It follows from (2) that $R(L(u_j))=R(Au_{j+1})=R(u_{j+1})=J_j$.
By Lemma \ref{2} (1), we get an isomorphism  $\varphi_{u\!_j, _1}:
\Hom_A(I_j, A)\simeq J_j$ of $k$-vector spaces.  In fact, we can
check directly that $\varphi_{u\!_j, _1}$ is an isomorphism of
$A^{op}$-modules. This proves (3).

(4) Note that $\Hom_A(I_j, I_i)= \Hom_A(Au_j, Au_i)\simeq u_jA\cap
Au_i = J_j\cap I_i$. To prove (4), there are three cases to be
considered.

Case 1: $j = i$.  By (1) and (2),  we conclude that $<x_2+q^jx_1,
x_1x_2, x_0x_1x_2>\subseteq I_j\cap J_j$. Since dim$_k(I_j)=4$ and
$x_2x_0-q^{j+1}x_0x_1\not\in I_j$, we get  dim$_k(I_j\cap J_j)=3$.
As a result, $ I_j\cap J_j=<x_2+q^jx_1, x_1x_2, x_0x_1x_2>$.

Case 2: $j=i-2$. Note that
$x_2x_0-q^{j+1}x_0x_1=x_2x_0-q^{i-1}x_0x_1$. But ${x_2+q^jx_1}\notin
 I_i$. It follows that  $I_i\cap J_j=<x_2x_0-q^{j+1}x_0x_1, x_1x_2,
x_0x_1x_2>$.

Case 3: $j\not\in\{i, i-2\}$. We claim that $I_i\cap J_j=<x_1x_2,
x_0x_1x_2>$. Obviously, $<x_1x_2,\, x_0x_1x_2>$ is contained in
$I_i\cap J_j$. Conversely, if $\lambda\in I_i\cap J_j$, then there
are elements $a_1, a_{20}, a_{21}, a_3, b_1, b_{20}, b_{21}
\,\mbox{and \,} b_3\in k$, such that $\lambda = a_1(x_2+q^jx_1)+
a_{20}(x_2x_0-q^{j+1}x_0x_1)+ a_{21}x_1x_2+ a_3x_0x_1x_2 =
b_1(x_2+q^ix_1)+ b_{20}(x_2x_0-q^{i-1}x_0x_1)+ b_{21}x_1x_2+
b_3x_0x_1x_2$. This implies that  $a_1=b_1, \,a_3 = b_3,\, a_{20} =
b_{20}, \,a_{21} = b_{21},\, a_1q^j = b_1q^i,$ and $ a_{20}q^{j+1} =
b_{20}q^{i-1} $. Consequently, $a_1 = a_{20} = 0$, which means that
$\lambda\in <x_1x_2, x_0x_1x_2>$. Thus $I_i\cap J_j=<x_1x_2,
x_0x_1x_2>$.

(5) The exact sequence $0\ra I_{j+1}\ra A\ra I_j\ra 0$ of
$A$-modules induces the following exact sequence of $k$-modules:
$$0\lra \Hom_A(I_j, I_i)\lra \Hom_A(A,
I_i)\lra \Hom_A(I_{j+1}, I_i)\lra \Ext_A^1(I_j, I_i)\lra 0.$$ By
(4), we have $$\dk\Hom_A(I_j, I_i)=
\left\{\begin{array}{ll} 3 & \mbox{if }  i\in\{j, j+2\},\\
2 & \mbox{otherwise}.\end{array} \right.$$ Since  $\dk(I_i)= 4$, we
have  $$\dk\Ext_A^1(I_j, I_i)=
\left\{\begin{array}{ll} 1 & \mbox{if } j\leq i\leq j+3,\\
0 & \mbox{otherwise}.\end{array} \right.$$ This proves (5).

$(6)$ By definition, we know that $A(1,u_i)=A$, and
$A(u_i,1)=R(u_{i+1})=J_i$.

$(7)$ It follows from $(4)$ and Lemma \ref{2}(2) that
$$
\dk A({u_j,u_i})=
\left\{\begin{array}{ll} 7 & \mbox{if } j=\{i-2, i\},\\
6 & \mbox{otherwise}.\end{array} \right.
$$
By definition, we know that $A({u_j,u_i})=\{a\in A~|~u_{j+1}au_i =
0\}$. It is not hard to see that
$$<x_1,x_2,x_0x_1,x_1x_2,x_2x_0,x_0x_1x_2>\subseteq A({u_j,u_i}).$$
Hence, if $j\not\in \{i-2, i\}$, then
$A({u_j,u_i})=<x_1,x_2,x_0x_1,x_1x_2,x_2x_0,x_0x_1x_2>$. If $j=i$,
then $u_{j+1}u_j=0$, and therefore $1\in A({u_j,u_j})$. Thus
$A({u_j,u_j})=<1,x_1,x_2,x_0x_1,x_1x_2,x_2x_0,x_0x_1x_2>$. If
$j=i-2$, then we can check that $u_{j+1}x_0u_{j+2}=0$. Thus, $x_0\in
A_{j,j+2}$. This shows that
$A({u_j,u_j+2})=<x_0,x_1,x_2,x_0x_1,x_1x_2,x_2x_0,x_0x_1x_2>$.
 $\square$

For higher cohomological groups, we have the following estimation.

\begin{Lem}\label{4}
Let  $t$ be an integer and $j$ a positive integer. Then

\smallskip
$(1)$ $\dk\Ext_A^j(I_0, I_t)=
\left\{\begin{array}{ll} 1 & \mbox{if } -1\leq t-j\leq 2,\\
0 & \mbox{otherwise}.\end{array} \right.$

\smallskip
$(2)$ $\dk\Ext_A^j(I_t, I_0)=
\left\{\begin{array}{ll} 1 & \mbox{if } -2\leq t+j\leq 1,\\
0 & \mbox{otherwise}.\end{array} \right.$

\smallskip
$(3)$ $\Ext_A^j(I_0, I_0)=0 \mbox{ for } j> 1.$
\end{Lem}

{\bf{Proof}}.  By Lemma \ref{lem2'}, we have  $\Ext_A^j(I_0,
I_t)\simeq \Ext_A^1(I_0, \Omega_A^{-j+1}(I_t))\simeq \Ext_A^1(I_0,
I_{t-j+1})$.  Now (1) follows from Lemma \ref{3}(5). Similarly, we
can prove (2). Clearly, (3) follows from (1) and (2). $\square$

\medskip
Here and subsequently,  $\delta_j$ stands for the canonical exact
sequence $ 0\ra I_{j+1}\ra A\ra I_j\ra 0$ in $A\modcat$ for each
$j\in \mathbb{Z}$.

\begin{Lem}\label{prop1}
Let $l\in \mathbb{Z}$ and $n\in \mathbb{N}$. Then
$$
\{j\in \mathbb{Z} \mid  \delta_j ~\mbox{is an\, $\add(A\sz
I_l)$-split sequence in A\modcat}\}= \{j\in\mathbb{Z}~|~ j>l+2
~\mbox{or}~ j<l-3\}.
$$
In particular, we have
$$
\{j\in \mathbb{Z}\mid \delta_j~ \mbox{is an
$\add(A\sz\bz\limits_{i=0}^{n}I_i)$-split sequence in $A\modcat$}\}=
\{j\in\mathbb{Z}~|~ j>n+2 ~\mbox{or}~ j<-3\}.
$$
\end{Lem}

{\bf Proof}. For any $ j\in \mathbb{Z}$, we know that $\delta_j$ is
an $\add(A\sz I_l)$-split sequence in $A\modcat$ if and only if
$\Ext_A^1(I_l, I_{j+1})=\Ext_A^1(I_j, I_l)=0$, which is equivalent
to the condition that $ j+1\not\in [l, l+3] \mbox{ and } j\not\in
[l-3, l]$ by Lemma \ref{3}(5). Thus we have (1). Clearly, (2)
follows from (1) immediately. $\square$

\medskip
The following  result can be directly deduced from the work of Hu
and Xi in \cite{hx2, HAU}.

\begin{Lem}\label{lem5}
Let $B$ be a $k$-algebra. Let   $Y$ and $M$ be $B$-modules with $M$
a generator for $B\modcat$. If $\Ext_B^1(M, \Omega_B(Y))=
\Ext_B^1(Y, M)=0$, then the endomorphism algebras $\End_B(M\sz Y)$
and $\End_B(M\sz \Omega_B(Y))$  are derived equivalent.  If, in
addition, $\Ext_B^2(M, \Omega_B(Y))= \Ext_B^2(Y, M)=0$, then the
$\{0, 1\}$-Auslander-Yoneda algebras $ \E_B^{\{0,1\}}(M\sz Y)$ and
$\E_B^{\{0,1\}}(M\sz \Omega_B(Y))$ are derived equivalent.

\end{Lem}

Having made the previous preparations, now we can prove Theorem
\ref{th1}.

{\bf Proof of Theorem \ref{th1}}. Let $m\geq n+4$. Set $M :=\,A\sz
I$ with $I=\bz\limits_{i=0}^{n}I_i$, and $V_m: = M \oplus I_m$.

(1) By Lemma \ref{3}(5), we know that  $\Ext_A^1(M, I_m)=
\Ext_A^1(I_m, M)= 0$. Clearly, we have
$$\dk(\Lambda_m^{\{0\}})=\dk\End_A(M)+\dk\Hom_A(M,I_m)+\dk\Hom_A(I_m,
M)+\dk\End_A(I_m) $$and
$$\dk(\Lambda_m^{\{0,1\}})= \dk(\Lambda_m^{\{0\}})+ \dk\Ext_A^1(M, M)+
\dk\Ext_A^1(I_m, I_m).$$ By Lemma \ref{3},   we get
$$\dk\End_A(I_m)=3,\; \dk\Ext_A^1(I_m, I_m)=1,
\;\dk\Hom_A(M,I_m)= \dk\Hom_A(I_m, M)=2n+6.$$  It follows that
$\dk{(\Lambda_m^{\Phi})}= \dk{(\Lambda_{m+1}^{\Phi})}$.

\smallskip
(2) We first show that  $\gd(\Lambda_m^{\{0\}})=\infty$. By Lemma
\ref{3}\,(5),  we have $\Ext_A^1(V_m, I_j)=0 $ for any $j<0$. Note
that, for any $t<j<0$, there is a long exact sequence
$$
0\lra I_j\lra A\lra A\lra \cdots\lra A\lra I_t\lra 0.
$$
It follows that the induced sequence
$$
0\lra \Hom_A(V_m, I_j)\lra \Hom_A(V_m, A)\lra \cdots\lra \Hom_A(V_m,
A)\lra \Hom_A(V_m, I_t)\lra 0
$$ is exact. Since
$\Hom_A(V_m, A)$ is a projective-injective indecomposable
$\Lambda_m^{\{0\}}$-module, we have
$\id_{\Lambda_m^{\{0\}}}\Hom_A(V_m, I_j) = \infty$ for all $j<0$.
Hence $\gd(\Lambda_m^{\{0\}})=\infty$. Note that there is a
canonical surjective homomorphism $\pi: \Lambda_m^{\{0,1\}}\ra
\Lambda_m^{\{0\}}$ of algebras. Thus every
$\Lambda_m^{\{0\}}$-module can be regarded as a
$\Lambda_m^{\{0,1\}}$-module. In addition, $\E_A^{\{0,1\}}(V_m, A)=
\Hom_A(V_m, A)$. It follows that
$\id_{\Lambda_m^{\{0,1\}}}\Hom_A(V_m, I_j) = \infty$ for all $j<0$.
This yields $\gd(\Lambda_m^{\{0,1\}})=\infty$.

\smallskip
(3) Recall a classical result on dominant dimension: Let $B$ an
algebra and $Y$ be a generator-cogenerator for $B\modcat$. Suppose
that $s$ is a non-negative integer.  Then $\dm(\End_B(Y))=s+2$ if
and only if $\Ext_B^i(Y, Y)=0$ for all $i$ with $1\leq i\leq s$, but
$ \Ext_B^{s+1}(Y, Y)\neq 0$. In our case, we take $Y:=V_m$ and
$s=0$. By Lemma \ref{3}(5), we know that $\Ext_A^1(I_0, I_0)\neq 0$,
which means that $\Ext_A^1(V_m, V_m)\neq 0$.  Note that $V_m$ is a
generator-cogenerator for $A\modcat$. Thus $\dm(\Lambda_m^{\{0\}})=
2$. By Lemma \ref{lem 2.13}, we have $\dm(\Lambda_m^{\{0,1\}})= 0$.

(4) Consider the exact sequence
$$
\delta_m: 0\lra I_{m+1}\lra A\lra I_m\lra 0$$ in $ A\modcat$. Since
$m\geq n+4$, it follows from Lemma \ref{3}(5) and Lemma
\ref{lem2'}(4) that $\Ext_A^1(M, I_{m+1})= \Ext_A^1(I_{m+1}, M)=
\Ext_A^1(I_m, M)= \Ext_A^1(M, I_{m})= 0$. Note that $A$ is
self-injective. By Lemma \ref{lem5}, we conclude that the algebras
$\Lambda_m^{\Phi}$ and $\Lambda_{m+1}^{\Phi}$ are derived-equivalent
for $\Phi=\{0\}$ or $\{0, 1\}$.

(5) It follows from Lemma \ref{lem2'} that  $\Omega_A(I_j)=I_{j+1}$
for each $j\in \mathbb{Z}$ and that the $A$-modules $I_j$  are
pairwise non-isomorphic for all $j\in \mathbb{Z}$. Now, we define
$W:= I_0$ and $W_n:=\oplus_{0\le j\le n} I_j$. Then, by Corollary
\ref{coro15}, the algebras $\Lambda_l^{\Phi}$ and $\Lambda_m^{\Phi}$
are not stably equivalent of Morita type if $l>m$. Thus the proof is
completed. $\square$

\medskip
In the rest of this section, we consider the special case: $n=0$ and
$\Phi=0$ in Theorem \ref{th1}.  For convenience, we set $\Lambda_m
:= \End_A(A\oplus I_0\oplus I_m)$ for $m\in\mathbb{Z}$, and define
$C:=<1,x_1,x_2,x_0x_1,x_1x_2,x_2x_0,x_0x_1x_2>,$
$T:=<x_1,x_2,x_0x_1,x_1x_2,x_2x_0,x_0x_1x_2>,$ and $S:=T\oplus <x_0>
$. Note that they all are subspaces of $A$.

\begin{Prop} Let $m$ be an  integer. Then

$(1)$ If $m\neq 2$, then $\Lambda_m$ is isomorphic to the algebra
$$
{\rm{M}}_{A}(1,u_0,u_m):= \left(
   \begin{array}{ccc}
     A   & A/I_1 & A/I_{m+1} \\
     J_0 & C/I_1 & T/I_{m+1} \\
     J_m & T/I_1 & C/I_{m+1} \\
   \end{array}
 \right).
$$

$(2)$ $\Lambda_2$ is isomorphic to the algebra
$$
 {\rm{M}}_{A}(1,u_0,u_2):=\left(
   \begin{array}{ccc}
     A   & A/I_1 & A/I_{3} \\
     J_0 & C/I_1 & S/I_{3} \\
     J_2 & T/I_1 & C/I_{3} \\
   \end{array}
 \right).$$

$(3)$ Suppose $m\geq 3$. Then, for any $l>m$, the algebras
$\Lambda_l$ and $\Lambda_m$ are derived-equivalent, but not stably
equivalent of Morita type. \label{prop3}
\end{Prop}

{\bf Proof}. (1) and (2) follow easily from Lemma \ref{2} and Lemma
\ref{3}, while $(3)$  can be concluded from Lemma \ref{prop1}, Lemma
\ref{lem5} and Corollary \ref{coro15}. $\square$

\bigskip
For each positive integer $m\geq 3$, the algebra $\Lambda_m$ is
given by the following quiver $Q$ with relations $\rho_m$:
$$
Q:\quad \xymatrix{\bullet\ar@<0.4ex>[r]^{\gamma_{_m}}
&\bullet\ar@<0.4ex>[l]^{\delta}^(1){3}^(0){1}\ar@(ul,ur)^{\alpha}\ar@<0.4ex>[r]^{\beta}
&\bullet\ar@<0.4ex>[l]^{\gamma_{_0}}^(0){2}}
$$

$$\begin{array}{ll} \rho_m: &
\alpha^2=\gamma_{_0}\beta\gamma_{_0}\alpha\beta
=\gamma_{_m}\alpha\delta\gamma_{_m}\delta=0;\\
& \\
  & \left(
   \begin{array}{ccc}
     \alpha  \beta  \gamma_{_0} \\
     \alpha  \delta  \gamma_{_m} \\
   \end{array}
 \right)=\frac{1}{q-q^{m+1}}
 \left(
   \begin{array}{cc}
     q^{m+2}-1 & 1-q^2 \\
     q^{m+2}-q^m & q^m-q^2 \\
   \end{array}
 \right)
 \left(
   \begin{array}{c}
     \beta\gamma_{_0}\alpha \\
     \delta\gamma_{_m}\alpha \\
   \end{array}
 \right); \\
 & \\
 & \frac{\beta\gamma_{_0}\beta}{1-q}=\frac{\delta\gamma_{_m}\beta}{q^m-q},\;~
\frac{\beta\gamma_{_0}\delta}{1-q^{m+1}}=\frac{\delta\gamma_{_m}\delta}{q^m-q^{m+1}};\quad
\frac{\gamma_{_0}\beta\gamma_{_0}}{1-q}=\frac{\gamma_{_0}\delta\gamma_{_m}}{1-q^{m+1}},\;~
\frac{\gamma_{_m}\beta\gamma_{_0}}{q^m-q}=\frac{\gamma_{_m}\delta\gamma_{_m}}{q^m-q^{m+1}}.
\end{array} $$

\bigskip The Cartan matrix of $\Lambda_m$ for $m\ge 3$ is
$$ C= \left(
   \begin{array}{ccc}
     8& 4 & 4 \\
     4& 3 & 2\\
     4& 2& 3 \\
   \end{array}
 \right),
$$ which is symmetric. Moreover, there is an anti-automorphism on
$\Lambda_m$ for $(m\ge 3)$, which is given by $$e_1\mapsto
e_1,\,e_2\mapsto e_3,\,e_3\mapsto e_2,\,\beta\mapsto\gamma_m,\,
\gamma_m\mapsto q^m\beta,\,\alpha\mapsto\alpha,\,
\delta\mapsto\gamma_0,\, \gamma_0\mapsto\delta.$$

It follows from Proposition \ref{prop3} that $\Lambda_t$, $t\ge 3$,
are pairwise derived-equivalent, but not stably equivalent of Morita
type.

Note that the Cartan matrix of $\Lambda_2$ is not symmetric. Thus
$\Lambda_2$ is not derived-equivalent to $\Lambda_m$ for $m\ge 3$
since the Cartan matrices of two derived equivalent algebras are
congruent over $\mathbb Z$, and therefore derived equivalences
preserve the symmetry of Cartan matrices. We don't know whether
$\Lambda_1$ and $\Lambda_3$ are derived-equivalent or not.

It would be interesting to show that the family of algebras in
Theorem \ref{th1} or in Proposition \ref{prop3} are pairwise not
stably equivalent.

\bigskip
%{\bf Acknowledgements.} The research work of the corresponding
%author C.C.Xi is partially supported by the Fundamental Research
%Funds for the Central Universities (2009SD-17). The revised version
%of the first draft of this paper was partially done when C.C.Xi
%visited the Chern Institute of Mathematics in July, 2010, he would
%like to thank Professor Chengming Bai for his warm hospitality.

\bigskip
June 30, 2010; Revised: November 10, 2010
\end{document}